\documentclass{amsart} 

\usepackage{amssymb} 
\usepackage{amsmath} 
\usepackage{amscd} 
\usepackage{latexsym} 
\usepackage{amsthm} 
\usepackage{verbatim} 

\let\cal\mathcal 
\def\Ascr{{\cal A}} 
 
\def\Cscr{{\cal C}} 
\def\Dscr{{\cal D}} 
\def\Escr{{\cal E}} 
\def\Fscr{{\cal F}} 
 
\def\Hscr{{\cal H}} 
\def\Iscr{{\cal I}} 
\def\Jscr{{\cal J}} 
 
\def\Lscr{{\cal L}} 
\def\Mscr{{\cal M}} 
\def\Nscr{{\cal N}} 
\def\Oscr{{\cal O}} 
\def\Pscr{{\cal P}} 
 
\def\Rscr{{\cal R}} 
\def\Sscr{{\cal S}} 
 
\def\Uscr{{\cal U}}

\def\Xscr{{\cal X}} 
\def\Yscr{{\cal Y}} 
 
\let\blb\mathbb 
\def\CC{{\blb C}}

\def \PP{{\blb P}} 
 
\def \ZZ{{\blb Z}} 
 
\def \NN{{\blb N}}

\def\id{\text{id}} 
\def\Id{\operatorname{id}} 
\def\pr{\mathop{\text{pr}}\nolimits}

\def\Res{\operatorname{Res}} 
 
\def\Lotimes{\overset{L}{\otimes}} 
\def\opp{{\operatorname{opp}}} 
\def\quot{/\!\!/} 
 
\def\Mod{\operatorname{Mod}} 
\def\mod{\operatorname{mod}} 
 
\def\fpd{\operatorname{fpd}} 
\def\fid{\operatorname{fid}}

\def\Lie{\mathop{\text{Lie}}}

\def\length{\mathop{\text{length}}}

\def\coh{\mathop{\text{\upshape{coh}}}}

\def\Rep{\operatorname {Rep}} 
 
\def\Gl{\operatorname {Gl}}

\def\Ext{\operatorname {Ext}} 
\def\Hom{\operatorname {Hom}} 
 
\def\RHom{\operatorname {RHom}} 
\def\uRHom{\operatorname {R\mathcal{H}\mathit{om}}} 
 
\def\cd{\operatorname {cd}}

\def\coker{\operatorname {coker}} 
\def\ker{\operatorname {ker}} 
 
\def\Tor{\operatorname {Tor}}

\def\id{{\operatorname {id}}}

\def\rk{\operatorname {rk}}

\def\Pic{\operatorname {Pic}} 
\def\gkdim{\operatorname {GK dim}} 
 
\def\gldim{\operatorname {gl\,dim}} 
 
\def\r{\rightarrow}

\DeclareMathOperator{\Proj}{Proj} 
 
\DeclareMathOperator{\Ind}{Ind}

\DeclareMathOperator{\Tors}{Tors} 
\DeclareMathOperator{\tors}{tors}

\DeclareMathOperator{\Aut}{Aut}

\newtheorem{lemma}{Lemma}[section] 
\newtheorem{proposition}[lemma]{Proposition} 
\newtheorem{theorem}[lemma]{Theorem} 
\newtheorem{corollary}[lemma]{Corollary} 
 
\newtheorem{convention}[lemma]{Convention}

\newtheorem{lemmas}{Lemma}[subsection] 
\newtheorem{propositions}[lemmas]{Proposition} 
\newtheorem{theorems}[lemmas]{Theorem} 
\newtheorem{corollarys}[lemmas]{Corollary}

\theoremstyle{definition}

\newtheorem{definitions}[lemmas]{Definition}

{ 
 
\newtheorem{step}{Step}

} 

\theoremstyle{remark} 

\newtheorem{remark}[lemma]{Remark} 
\newtheorem{remarks}[lemmas]{Remark}

\newdimen\uboxsep \uboxsep=1ex 
\def\uboxn#1{\vtop to 0pt{\hrule height 0pt depth 0pt\vskip\uboxsep 
\hbox to 0pt{\hss #1\hss}\vss}} 

\def\uboxs#1{\vbox to 0pt{\vss\hbox to 0pt{\hss #1\hss} 
\vskip\uboxsep\hrule height 0pt depth 0pt}}

\numberwithin{equation}{section} 

\let\cal\mathcal 
\def\Ascr{{\cal A}} 
 
\def\Cscr{{\cal C}} 
\def\Dscr{{\cal D}} 
\def\Escr{{\cal E}} 
\def\Fscr{{\cal F}} 
 
\def\Hscr{{\cal H}} 
\def\Iscr{{\cal I}} 
\def\Jscr{{\cal J}} 
 
\def\Lscr{{\cal L}} 
\def\Mscr{{\cal M}} 
\def\Nscr{{\cal N}} 
\def\Oscr{{\cal O}} 
\def\Pscr{{\cal P}} 
 
\def\Rscr{{\cal R}} 
\def\Sscr{{\cal S}} 
 
\def\Uscr{{\cal U}}

\def\Xscr{{\cal X}} 
\def\Yscr{{\cal Y}}

\let\blb\mathbb 
\def\CC{{\blb C}}

\def\PP{{\blb P}} 
 
\def\ZZ{{\blb Z}} 
 
\def\NN{{\blb N}}

\def\K0{{K_{0}(\blb P^{2}_{q})}} 

\def\pr{\mathop{\text{pr}}\nolimits} 
\def\id{\operatorname {id}} 
\def\Lotimes{\overset{{\bf L}}{\otimes}}
\def\LotimesA{\overset{{\bf L}}{\otimes}_{A}}
\def\LotimesD{\overset{{\bf L}}{\otimes}_{D}}
\def\Mod{\operatorname{Mod}} 
\def\mod{\operatorname{mod}} 
\def\coh{\mathop{\text{\upshape{coh}}}} 
 
\def\Ext{\operatorname {Ext}} 
\def\Hom{\operatorname {Hom}} 
 
\def\RHom{\operatorname {{\bf R}Hom}} 
\def\cd{\operatorname {cd}}

\def\coker{\operatorname {coker}} 
\def\ker{\operatorname {ker}} 
\def\Tor{\operatorname {Tor}} 
\def\Pic{\operatorname {Pic}} 
\def\gldim{\operatorname {gl\,dim}} 
\def\GrMod{\operatorname {GrMod}} 
\def\grmod{\operatorname {grmod}} 
\def\Tails{\operatorname {Tails}} 
\def\tails{\operatorname {tails}} 
\def\Qcoh{\operatorname {Qcoh}} 
 
\def\Fract{\operatorname {Fract}} 
\def\sh{\operatorname{sh}} 
\def\Skly{\operatorname{Skl}} 
 
\def\P2q{\operatorname {{\blb P}^{2}_{q}}} 
\def\D{\operatorname{D}} 
\def\Res{\operatorname {Res}} 
\def\Ind{\operatorname {Ind}} 
\def\cone{\operatorname {cone}}

\def\r{\rightarrow}

\def\lr{\longrightarrow}

\newdimen\uboxsep \uboxsep=1ex 
\def\uboxn#1{\vtop to 0pt{\hrule height 0pt depth 0pt\vskip\uboxsep 
\hbox to 0pt{\hss #1\hss}\vss}} 

\def\uboxs#1{\vbox to 0pt{\vss\hbox to 0pt{\hss #1\hss} 
\vskip\uboxsep\hrule height 0pt depth 0pt}}


\author{Koen De Naeghel and Michel Van den Bergh} 

\address{Departement WNI\\Limburgs Universitair 
Centrum\\ Universitaire Campus\\ Building D\\ 3590 
Diepenbeek\\ Belgium} 
\thanks{The second author is a director of research at the FWO} 

\email[K. De Naeghel]{koen.denaeghel@luc.ac.be} 
\email[M. Van den Bergh]{michel.vandenbergh@luc.ac.be} 
\date{April 1, 2003} 
\title[Ideal classes of three dimensional Sklyanin algebras]{Ideal classes 
of three dimensional Sklyanin algebras} 

\begin{document} 

\begin{abstract} 
In this paper we classify graded reflexive ideals, up to isomorphism 
and shift, in certain three dimensional Artin-Schelter regular 
algebras. This classification is similar to the classification of right 
ideals in the first Weyl algebra, a problem that was completely 
settled recently. 
The situation we consider is substantially more complicated however. 
\end{abstract} 

\maketitle 
\tableofcontents 
\section{Introduction} 
This paper is motivated by the recent developments on the 
classification of right ideals in the first Weyl algebra. 
We start by recalling the (for now) definite result in this subject, as it was 
formulated  by Berest and Wilson. 
\begin{theorem} \cite{BW0} 
\label{ref-1.1-0} 
Let $A_1=\CC\langle x,y\rangle/(yx-xy-1)$ be the first Weyl algebra. 
Put $G=\Aut(A_1)$ and let $\Rscr$ be the set of isomorphism classes of 
right $A_1$-ideals. Then the orbits of the natural $G$-action on 
$\Rscr$ are indexed by $\NN$, and the orbit corresponding 
to $n\in \NN$ is in natural bijection with the $n$'th Calogero-Moser 
space 
\begin{equation} 
\label{ref-1.1-1} 
C_n=\{X,Y\in M_n(\CC)\mid \rk(YX-XY-\Id)=1\} /\Gl_n(\CC) 
\end{equation} 
where  $\Gl_n(\CC)$ acts  by simultaneous conjugation on $(X,Y)$. 
\end{theorem} 
The fact $\Rscr/G\cong \NN$ has also been proved by Kouakou in his 
(unpublished) PhD-thesis \cite{kouakou}. 

The first proof of Theorem \ref{ref-1.1-0} used the fact that there is a 
description of $\Rscr$  in terms of 
the adelic Grassmanian (due to Cannings and Holland 
\cite{CH1}). 
Using methods from integrable systems, Wilson established a relation between 
the adelic Grassmanian and the Calogero-Moser spaces \cite{wilson}. 

In \cite{BW} Berest and Wilson gave a new proof of Theorem \ref{ref-1.1-0} using 
noncommutative algebraic geometry  \cite{AZ,VW} 
(some of the original proofs in \cite{BW} were slightly simplified by 
the second author).  See also \cite{Ginzburg,Orlov1}. 
That an 
approach based on noncommutative geometry should be possible was in 
fact anticipated very early by Le 
Bruyn who in \cite{lebruyn1} already came very close to proving Theorem 
\ref{ref-1.1-0}. 

Let us briefly indicate how the methods of noncommutative algebraic 
geometry  may be used to prove 
Theorem \ref{ref-1.1-0}. We introduce the \emph{homogenized Weyl algebra} 
$H=k\langle x,y,z\rangle/(zx-xz,zy-yz,yx-xy-z^2)$ and then we consider 
$A_1$ as the coordinate ring of an open affine part of a 
noncommutative space $\P2q$, with ``homogeneous coordinate ring'' $H$ 
(see below for more precise definitions). The problem of 
describing $\Rscr$ then becomes equivalent to describing certain 
objects on the ``noncommutative projective plane'' $\P2q$. 
Objects on $\P2q$ have finite dimensional cohomology groups and these 
may be used to define moduli spaces, just as in the ordinary 
commutative case. 

The current paper starts from the observation that there are many more 
noncommutative projective planes than just the one associated to the 
Weyl algebra (this is in fact a fairly degenerate one) 
\cite{ATV1,ATV2,Bondal}. So below we let $\P2q$ be  a so-called 
``elliptic'' quantum projective plane. 
By definition $\P2q$ is noncommutative projective scheme which has 
as homogeneous 
coordinate ring a graded ring $A$ with generators  $x,y,z$ (in 
degree one) satisfying the relations 
\begin{equation} 
\left\{ 
\begin{array}{l} 
ayz + bzy + cx^{2} =0 \\ 
azx + bxz + cy^{2}=0\\ 
axy + byx + cz^{2}=0 
\end{array} \right. 
\end{equation} 
where $a,b,c$ are generic scalars (see 
below). 

Such a $\P2q$ is called an elliptic quantum plane because there is an 
inclusion (in a noncommutative geometry sense) $E\hookrightarrow 
\P2q$ where $E$ is a smooth (commutative) elliptic curve. 

We let $\Rscr$ be the set of reflexive graded $A$-ideals, 
considered up 
to isomorphism and shift of grading. We may think of the elements of 
$\Rscr$ as line bundles on $\P2q$. In this paper we prove the 
following result (see Theorem \ref{ref-5.5.5-65} below). 
\begin{theorem} \label{ref-1.2-2} There exist smooth affine 
varieties $D_n$ of dimension 
$2n$ such that $\Rscr$ is naturally in bijection with 
$\coprod_n D_n$. 
\end{theorem} 
We would like to think of the $D_n$ as elliptic Calogero-Moser 
spaces. We show below that $D_0$ is a point and $D_1$ is the 
complement of $E$ under a natural embedding in $\PP^2$. 
\begin{remark}
In fact $D_n$ is connected, which we will prove in a subsequent paper
\cite{DV2}. 
\end{remark}
\begin{remark}
A theorem similar to Theorem \ref{ref-1.2-2} has been announced  by \cite{NS}. 
They work in a more general setting where the associated automorphism $\sigma$
of $E$ may have finite order. 
\end{remark}
The reader will notice that Theorem \ref{ref-1.2-2} is weaker than Theorem 
\ref{ref-1.1-0} but this is probably unavoidable. Although we have a fairly 
succinct description of the varieties $D_n$ 
(see \eqref{ref-5.9-64} below) it is not as explicit as 
\eqref{ref-1.1-1}. And very likely  $D_n$ can 
also not be viewed in a natural way as the orbit of a group. 

Our proof of Theorem \ref{ref-1.2-2} is similar in spirit to the proof 
of Theorem \ref{ref-1.1-0}. However it is substantially more involved. The 
reason for this is that the proofs for the Weyl algebra rely heavily on the 
fact that  $H$ contains a central element in 
degree one (namely $z$) and the lowest central element 
in $A$ has degree three. 

We also have a result which explicitly describes the elements of 
$\Rscr$.  Recall that a line module over $A$ is a graded $A$-module of the 
form $A/uA$ with $u\in A_1-\{0\}$. 

The following theorem can be deduced easily from Theorem 
\ref{ref-5.6.6-72} below. 
\begin{theorem} Let $I\in \Rscr$.   Then there exists 
an $m\in\NN$ together with a monomorphism $I(-m)\hookrightarrow 
A$ such that there exists a filtration of reflexive graded $A$-ideals 
$A=M_0\supset M_1\supset \cdots 
\supset M_u=I(-m)$ with the property that the $M_i/M_{i+1}$ are shifted 
line modules, up to finite length modules. 
\end{theorem} 
It seems plausible that this result may be used to obtain an analogue 
of the Cannings-Holland classification of ideals in the Weyl algebra 
(see \cite{CH1}) 
but we have not sorted out the details. We hope to come back on this 
in a subsequent paper. 
\section{Preliminaries} 
\label{ref-2-3} 
Throughout we work over an algebraically closed field $k$ of 
characteristic zero. In 
this section we recall some basic notions of noncommutative 
projective geometry. These are collected from 
\cite{AZ,MS,smith1,VdBSt,StICM,VdB19}. We use 
the following convention: 
\begin{convention} 
\label{ref-2.1-4} 
If $\mathrm{XyUvw}(\cdots)$  denotes an 
abelian category then 
$\mathrm{xyuvw}(\cdots)$ denotes the 
full subcategory  of $\mathrm{XyUvw}(\cdots)$ consisting of noetherian 
objects. 
\end{convention} 
To simplify the notations we often use implicitly the following result 
\begin{lemma} 
\label{ref-2.2-5} 
Assume that $\Cscr$ is a locally noetherian category and 
$\Cscr_f$ is the full subcategory of $\Cscr$ consisting of 
noetherian objects. Then the natural map $D^b(\Cscr_f)\r 
D^b_{\Cscr_f}(\Cscr)$ is an equivalence of categories. 
\end{lemma} 
\begin{proof} This follows for example from the dual of \cite[1.7.11]{KSI}. 
\end{proof} 
\subsection{Graded algebras and modules} 
Let $A=\oplus_{i\in \ZZ} A_i$ be a $\ZZ$-graded algebra. If $A_{i} = 0$ 
for all $i < 0$ we say that $A$ is {\it positively graded}, and if in addition 
$A_{0} = k$ we say that $A$ is {\it connected}. Any graded connected Noetherian 
$k$-algebra $A$ is {\it locally finite}, 
i.e. $\dim_{k}A_{i} < \infty$, for all $i \in \ZZ$.\\ 
We write $\GrMod(A)$ for the category of graded right 
$A$-modules with morphisms the $A$-module homomorphisms of degree 
zero. 
Let $M$ be a graded right $A$-module. We use the notation (for all $n 
\in \ZZ$) 
$M_{\geq n} = \oplus_{d \geq n}^{}{M_{d}}$ and 
$M_{\leq n} = \oplus_{d \leq n}^{}{M_{d}}$. 

We say that $M$ is {\it left} (resp. {\it right}) {\it bounded} if $M _{\leq 
n} = 0$ (resp. $M_{\geq n} = 0$) for some $n \in \ZZ$. 
For any integer $n$, define $M(n)$ 
as the graded $A$-module that is equal to $M$ with its original $A$ 
action, but which is graded by 
$M(n)_{i} = M_{n + i}$. 
We refer to the 
functor $M \mapsto M(n)$ as 
the $n$-th {\it shift functor}. 

Since $\GrMod(A)$ 
is an abelian  category with enough injective 
objects we may define the functors $\Ext_{A}^{n}(M,-)$ on 
$\GrMod(A)$ as the right derived functors of $\Hom_{A}(M,-)$. It 
is convenient to write (for $n \geq 0$) 
\[ 
\underline{\Ext}_{A}^{n}(M,N) := \bigoplus_{d \in 
\ZZ}^{}{\Ext_{A}^{n}(M,N(d))}; 
\] 
whence $\underline{\Ext}_{A}^{n}(M,-)$ are the right derived functors of 
$\underline{\Ext}_{A}^{0}(M,-) := \underline{\Hom}_{A}(M,-)$, for $n 
\geq 1$.

Finally, recall that a module $M \in \GrMod(A)$ is {\it reflexive} 
if $M^{**} = M$ where $M^* = \underline{\Hom}_{A}(M,A)$ is the graded 
dual of $M$. 
\subsection{Tails} 
Let $A$ be a Noetherian connected graded $k$-algebra. We denote by 
$\tau$ the functor that sends a graded right $A$-module to the 
the sum of all its finite dimensional 
submodules. 

Denote by $\Tors(A)$ the full subcategory of $\GrMod(A)$ 
consisting of all modules such that $\tau M=M$  and write $\Tails(A)$ for the 
quotient category \linebreak $\GrMod(A)/\Tors(A)$. 
We write $\pi : \GrMod(A) \r \Tails(A)$ for the (exact) quotient 
functor. By localization theory \cite{stenstrom} $\pi$ has a right 
adjoint which we 
denote by $\omega$. It is well-known that $\pi\circ\omega=\id$. 
The object $\pi A$ in $\Tails(A)$ will be 
denoted by $\Oscr$ and it is again easy to see that 
$\omega=\underline{\Hom}_{\Tails}(\Oscr,-)$. 
The objects in $\Tails(A)$ will be denoted by script letters, like $\Mscr$. 

The shift functor induces an 
automorphism $\sh : \Mscr \mapsto \Mscr(1)$ on $\Tails(A)$ which we 
also call the shift functor (in analogy with algebraic geometry it 
should perhaps be called the ``twist'' functor). 

When there is no possible confusion we write $\Hom$ 
instead of $\Hom_{A}$ and \linebreak
$\Hom_{\Tails(A)}$. The context will make clear 
in which category we work. 

If  $\Mscr \in \Tails(A)$ then $\Hom(\Mscr, -)$ is left 
exact, so we may define its right derived functors $\Ext^{n}(\Mscr, -)$. 
We also use the notation 
\[ 
\underline{\Ext}^{n}(\Mscr,\Nscr) := \bigoplus_{d \in 
\ZZ}^{}{\Ext^{n}(\Mscr,\Nscr(d))} 
\] 
and we set $\underline{\Hom}(\Mscr,\Nscr) = 
\underline{\Ext}^{0}(\Mscr,\Nscr)$. 

Our Convention \ref{ref-2.1-4} fixes the meaning of $\grmod(A)$, $\tors(A)$ 
and $\tails(A)$. 
It is easy to see that $\tors(A)$ consists of the 
finite dimensional graded $A$-modules. Furthermore 
$\tails(A)=\grmod(A)/\tors(A)$. 

If $M$ is finitely generated and $N$ is arbitrary we have 
\begin{equation} 
\label{ref-2.1-6} 
\underline{\Ext}^{n}(\pi M,\pi N) \cong \lim_{\lr} 
\underline{\Ext}^{n}_{A}(M_{\geq m},N). 
\end{equation} 
If $M$ and $N$ are 
both finitely generated, then (\ref{ref-2.1-6}) implies 
\[ 
\pi M \cong \pi N \mbox{ in } \tails(A) \hspace{0.25cm} \Leftrightarrow 
\hspace{0.25cm} M_{\geq n} \cong 
N_{\geq n} 
\mbox{ in $\grmod(A)$ for some $n$} 
\] 
explaining the word ``tails''. 

For $M \in \GrMod(A)$ there is an exact 
sequence (see \cite{AZ}, Proposition 7.2) 
\begin{eqnarray} \label{ref-2.2-7} 
0 \r \tau M \r M \r \omega \pi M \r \lim_{\lr} 
\underline{\Ext}_{A}^{1}(A / A_{\geq n},M) \r 0. 
\end{eqnarray} 
An object $\Mscr \in \Tails(A)$ is said to be {\it reflexive} if $\Mscr = \pi M$ 
for some reflexive $M \in \GrMod(A)$.

We say that $A$ 
satisfies condition $\chi$ if $\dim_{k}\Ext^{j}(k,M) < \infty$ for all $j$ 
and all $M \in \grmod(A)$. \\ 
In case $A$ satisfies condition $\chi$ then for every $M \in \grmod(A)$ 
the cokernel of the map $M \r \omega \pi M$ in the exact sequence 
(\ref{ref-2.2-7}) is right bounded. In particular, for $M\in \grmod(A)$ 
we have $M_{\geq d} \cong (\omega \pi M)_{\geq d}$ for some $d$. \\ 

Every graded quotient of a 
polynomial ring satisfies condition $\chi$ and so do most 
noncommutative algebras 
of importance. The condition is essential to get a theory for 
noncommutative schemes which resembles the commutative theory. 
\begin{propositions} \cite{AZ} Let $A$ be a right Noetherian 
connected $k$-algebra satisfying condition $\chi$. Then 
$\Ext^{j}(\Mscr,\Nscr)$ is finite dimensional for all $j$ and all \linebreak
$\Mscr,\Nscr \in \tails(A)$. 
\end{propositions}

\subsection{Serre duality} 
It was shown in \cite{YZ} that under reasonable hypotheses the 
category $\tails(A)$ satisfies a classical form of Serre duality. However we 
will need a 
stronger form of Serre duality introduced by Bondal and Kapranov in 
\cite{Bondal4}. Let $\Ascr$ be a $k$-linear $\Ext$-finite triangulated 
category. By this we mean that for all 
$\Mscr,\Nscr\in\Ascr$ we have $\sum_n\dim_k \Hom(\Mscr,\Nscr[n])<\infty$. 
The category $\Ascr$ is said to satisfy 
Bondal-Kapranov-Serre (BKS) duality  if there is an auto-equivalence 
$F:\Ascr\r \Ascr$ together with for all $A,B\in \Ascr$ natural 
isomorphisms 
\[ 
\Hom(A,B)\r \Hom(B,FA)' 
\] 
(where $(-)'$ denotes the $k$-dual). 

Let $\Cscr$ be an  abelian category. We say that $\Cscr$ has \emph{finite 
global dimension} if there exists an $n$ such that 
$\Ext^i_{\Cscr}(A,B)=0$ for all $A,B\in\Cscr$ and for all $i>n$. The 
minimal such $n$ is called the \emph{global dimension} of $\Cscr$. 

In this section we assume that $A$ is a connected graded noetherian 
ring over a  $k$. By $(-)'$ we denote the functor on graded 
vector spaces which sends $M$ to $\oplus_n M_{-n}^\ast$. If we use 
notations which refer to the left structure of $A$ then we adorn them 
with a superscript ``$\circ$''. 

We make the following additional assumptions on $A$ 
\begin{enumerate} 
\item $A$ satisfies $\chi$ and the functor $\tau$ has finite 
cohomological dimension. 
\item $A$ satisfies $\chi^\circ$ and the functor $\tau^\circ$ has finite 
cohomological dimension. 
\item $\tails(A)$ has finite global dimension. 
\end{enumerate} 
Note that if $A$ has finite global dimension then so does $\tails(A)$ 
by \eqref{ref-2.1-6}. 

Put $R=R\tau(A)'$. According to \cite{VdB16} $R$ is a complex of 
bimodules with finitely generated cohomology on the left and on the 
right, which in addition has finite injective dimension, also on the 
left and on the right. We now have the following result 
\begin{theorems} 
\label{ref-2.3.1-8} 
(Serre duality) The category $D^b(\tails(A))$ 
satisfies BKS-duality with Serre functor defined by 
\[ 
F(\pi M)=\pi(M\Lotimes R)[-1] 
\] 
\end{theorems} 
This result is certainly not unexpected but as far as we know a 
written proof does not exist in the literature. We prove a more general 
version of Theorem \ref{ref-2.3.1-8} in Appendix \ref{ref-A-74}.

\subsection{Projective schemes} 
We  use the definition of $\Proj A$ (for a noncommutative 
algebra $A$) suggested by Artin and Zhang (see \cite{AZ}). 
Let $A$ be a Noetherian graded $k$-algebra. 
We define the (polarized) projective scheme $\Proj A$ of $A$ as the 
triple $(\tails(A),\Oscr,\sh)$. 
In what follows we shall refer to the objects of 
$\tails(A)$ (resp. $\Tails(A)$) as the {\it coherent} (resp. 
{\it quasicoherent}) sheaves on 
$X = \Proj A$, even when $A$ is not commutative, and we shall use the notation 
$\coh(X) := \tails(A)$, \linebreak $\Qcoh(X) := \Tails(A)$. By analogy we sometimes 
write $\Oscr_X=\Oscr=\pi A$. 

The following definitions agree with the classical ones for 
projective schemes. \\ 
If $\Mscr$ is be a quasicoherent sheaf on $X = \Proj A$, we define the 
{\it cohomology groups} of $\Mscr$ by 
\[ 
H^{n}(X,\Mscr) := \Ext^{n}(\Oscr_{X},\Mscr). 
\] 
We refer to the graded right $A$-modules 
\[ 
\underline{H}^{n}(X,\Mscr) := \bigoplus_{d \in \ZZ}^{}{H^{n}(X,\Mscr(d))} 
\] 
as the {\it full cohomology modules} 
of $\Mscr$. \\ Finally, we mention  the 
{\it cohomological dimension} of $X$ 
\[ 
\cd X := \max \{ n \in \NN \mid  H^{n}(X, -) \neq 0 \}. 
\] 
It is easy to prove that 
\[ 
\cd X=\max(0,\cd \tau-1). 
\] 

\subsection{The Grothendieck group, the Euler form and Hilbert series} 
In this subsection $A$ will be a Noetherian connected graded $k$-algebra 
with finite global (homological) dimension. We recall some basic 
tools. 

\subsubsection{The Euler form} \label{ref-2.5.1-9} 
Let $\Cscr$ be an $\Ext$-finite $k$-linear abelian category of finite 
global dimension. 
We define 
\[ 
\chi(A,B)=\sum_i (-1)^i\dim \Ext^i_\Cscr(A,B) 
\] 
for $A,B\in \Cscr$. It is clear that $\chi$ defines a bilinear form 
$\chi:K_0(\Cscr)\times 
K_0(\Cscr)\r \ZZ$ which we call the \emph{Euler form} for $\Cscr$.

\subsubsection{Hilbert series} 
The Hilbert series of $M \in \grmod(A)$ is the Laurent power series 
\[ 
h_{M}(t) = \sum_{i = - \infty}^{+ \infty}{(\dim_{k}M_{i})t^{i}} \in\ZZ((t)). 
\] 
This definition makes sense since $A$ is right Noetherian. 

Let $M \in \grmod(A)$. Given a resolution 
\[ 
0 \r P^{r} \r \dots \r P^{1} \r P^{0} \r M \r 0 
\] 
we have 
\[ 
h_{M}(t) = \sum_{i=0}^{r}{(-1)^{i}h_{P^{i}}(t)}. 
\] 
Since $A$ is connected, left bounded graded right $A$-modules are 
projective if and only if they are free 
hence isomorphic to a sum of shifts of $A$. 
So if we write 
\[ 
P^{i} = \bigoplus_{j=0}^{r_{i}}{A(-l_{ij})} 
\] 
we obtain the formula 
\begin{eqnarray} \label{ref-2.3-10} 
q_{M}(t) = h_{M}(t) h_{A}(t)^{-1} 
\end{eqnarray} 
where $q_{M}(t)$ is the \emph{characteristic polynomial of $M$}, it is 
defined by 
\[ 
q_{M}(t) = \sum_{i=0}^{r}{(-1)^{i}} \sum_{j=0}^{r_{i}}{t^{l_{ij}}} 
\in \ZZ[t,t^{-1}]. 
\] 

\subsubsection{The Grothendieck group $K_{0}(X)$ and the rank 
function $K_{0}(X) \r \ZZ$}  
Set \linebreak $X = \Proj A$. 
If $\Mscr$ is a coherent sheaf on $X$, we denote by $[\Mscr]$ its image in 
$K_{0}(X)$. \\ 
The shift functor on $\coh(X)$ induces an  automorphism of $K_{0}(X)$. 
Following 
\cite{MS}, we view $K_{0}(X)$ as a $\ZZ[t,t^{-1}]$-module with $t$ 
acting as the shift functor \linebreak $\sh^{-1}: \Mscr \mapsto \Mscr(-1)$. 
Now $K_{0}(X)$ may be described in terms of the Hilbert series of $A$. 
\begin{theorems} [\cite{MS}, Theorem 2.3] \label{ref-2.5.1-11} 
Let $A$ be a Noetherian connected graded $k$-algebra of finite global 
dimension. 
Set $X = \Proj A$ and let $q = q_{k}(t)$. Then 
\[ 
K_{0}(X) \cong \ZZ[t,t^{-1}]\left/(q)\right. 
\] 
and for each $M \in \grmod(A)$, the isomorphism sends $[\pi M]$ to the 
characteristic polynomial $q_{M}(t)$ of $M$: 
\[ 
\theta:[\pi M] \mapsto \overline{q_{M}(t)} 
\] 
In particular, $[\Oscr(n)]$ is sent to $t^{-n}$. 
\end{theorems} 
Let us assume that $A$ is a domain generated in degree one, so it has a 
graded division ring of fractions 
(graded version of Goldie's Theorem; see \cite[Ch. C, Cor. I.1.7]{NVO1}) 
\[ 
\Fract(A) := \{ ab^{-1} \mid a,b \in A \mbox{ homogeneous, $b \neq 0$} \}. 
\] 
The degree zero component $k(X)$ of $\Fract(A)$ is a division algebra 
which is called the 
{ \it function field} of $X$. 
$\Fract(A)$ is isomorphic to a skew Laurent extention 
$k(X)[z,z^{-1};\sigma]$ where $z$ has degree one 
(see \cite{NVO1}, Chapter A, Corollary I.4.3). \\ 
The {\it rank} of a finitely graded right $A$-module $M$ is 
\[ 
\rk M = \dim_{k(X)}\left( M \otimes _{A} \Fract(A) \right)_{0} 
\] 
This also defines an additive rank fuction on $\coh(X)$ and hence a 
homomorphism $K_0(X)\r \ZZ$ also denoted by ``$\rk$''. Obviously 
$\rk \Oscr=1$ and $\rk \Mscr=\rk \Mscr(1)$.

\subsection{Artin-Schelter regular algebras} 
\label{ref-2.6-12} 
\begin{definitions} \cite{AZ} 
A connected graded $k$-algebra $A$ is called an {\it Artin-Schelter 
regular algebra of dimension} $d$ if it has the following properties: 
\begin{enumerate} 
\item[(i)] A has finite global dimension $d$; 
\item[(ii)] A has polynomial growth, that is, there exists positive real 
numbers $c, 
\delta$ 
such that $\dim_{k}A_{n} \leq cn^{\delta}$ for all positive integers $n$; 
\item[(iii)] A is Gorenstein, meaning there is an integer $l$ such that 
\[ 
\underline{\Ext}_{A}^{i}(k_{A},A) \cong 
\left \{ 
\begin{array}{ll} 
_{A}k(l) & \mbox{ if $i = d$,}\\ 
0 & \mbox{ otherwise.} 
\end{array} 
\right. 
\] 
where $l$ is called the {\it Gorenstein parameter} of $A$. 
\end{enumerate} 
\end{definitions} 
If $A$ is commutative, then the condition (i) already implies that $A$ is 
isomorphic 
to a polynomial ring $k[x_{1}, \dots x_{n}]$ with some positive 
grading. If in this case the grading is standard then $n = l$.\\ 
The Gorenstein property determines the full cohomology groups of 
$\Oscr$. 
\begin{theorems} \cite{AZ} \label{ref-2.6.2-13} 
Let $A$ be a Noetherian Artin-Schelter regular algebra of dimension 
$d = n + 1$, and let $X = \Proj A$. Then $\cd X = n$, and the full  
cohomology modules of $\Oscr = \pi A$ are given by 
\[ 
\underline{H}^{i}(X,\Oscr) \cong 
\left\{ 
\begin{array}{ll} 
A & \mbox{ if } i = 0 \\ 
0 & \mbox{ if } i \neq 0,n \\ 
A'(l) & \mbox{ if } i = n 
\end{array} 
\right. 
\] 
\end{theorems} 
Let $A$ be an Artin-Schelter regular algebra as in the previous 
Theorem and put $X=\Proj A$. It is easy to see that 
$A$ satisfies the hypotheses for Theoreom \ref{ref-2.3.1-8}.  In this 
case the Serre functor has a particularly simple form: 
indeed in \cite{AZ} it is shown that $R=(R^{n+1}\tau A)'\cong 
A[n+1](-l)$ as left $A$-modules and in \cite{VdB16}  it is proved that 
$R\tau A\cong R \tau^\circ A$ as complexes of bimodules. Thus we also have 
that $R=A[n+1](-l)$ as 
right $A$-modules. In other words $R=A_\phi[n+1](-l)$ where $\phi$ is some 
graded automorphism of $A$.  The automorphism $M\mapsto M_\phi$ of 
$\GrMod(A)$ passes to an automorphism $\Tails(A)$ for which we also 
use the notation $(-)_\phi$. 

We find the the following formula for the 
Serre functor on $\tails(A)$. 
\[ 
F\Mscr=\Mscr_\phi(-l)[n] 
\] 
From this we easily obtain: 
\begin{propositions} 
\label{ref-2.6.3-14} 
One has $\gldim \tails(A)= \gldim A-1$. 
\end{propositions} 
\begin{proof} As above put $\gldim A=n+1$. 
The inequality $\gldim \tails(A)\le n$ follows 
directly from BKS-duality and the above discussion.  The other 
inequality follows from Theorem \ref{ref-2.6.2-13}. 
\end{proof} 
\subsection{Dimension and multiplicity} 
\label{ref-2.7-15} 
Let $A$ be a noetherian Artin-Schelter regular algebra. 
If $0\neq M\in \grmod(A)$ then the 
\emph{Gelfand-Kirilov dimension} $\gkdim M$ of $M$ \cite{KL} can be computed 
as the order of the 
pole of $h_M(t)$ in $1$ \cite{ATV2} (in particular it is an 
integer). If $\gkdim M\le n$ then we define $e_n(M)$ as $\lim_{t\r 
1}(1-t)^n h_M(t)$. Clearly $e_n$ is additive on short exact 
sequences of objects with $\gkdim \le n$. We have $e_n(M)\ge 0$ and 
furthermore $e_n(M)=0$ if and only if $\gkdim M < n$. 
If $e_n(M)>0$ then we put 
$e(M)=e_n(M)$ and we call this the \emph{multiplicity} of $M$. 
If $u=\gkdim A$ and if $A$ is a domain generated in degree one then it 
is easy to see that 
$\rk M=e_u(M)/e_u(A)$. 

If $\Mscr=\pi M\neq 0$ then we put $\dim \Mscr=\gkdim M-1$ and 
$e_n(\Mscr)=e_{n+1}(M)$, $e(\Mscr)=e(M)$. 

An object in $\grmod(A)$ or $\tails(A)$ is said to be \emph{pure} if it 
contains no subobjects of strictly smaller dimension. It is \emph{critical} 
if all non-trivial subobjects have the same multiplicity. It is easy 
to see that if $M\in\grmod(A)$ is  pure or critical then so is $\pi M$, and 
conversely if $\Mscr\in \tails(A)$ is pure or critical then there 
exists a module $M\in\grmod(A)$ which has the corresponding property such that 
$\Mscr=\pi M$. 
\subsection{Three dimensional Artin-Schelter regular algebras} 
\label{ref-2.8-16} 
There exists a 
complete classification for Artin-Schelter regular algebras 
of dimension three \cite{ATV1,ATV2,Steph1,Steph2}: \label{ref-2.8-17} 
\begin{theorems} 
The Artin-Schelter regular algbras $A$ of dimension three can be 
classified. They are all Noetherian domains with Hilbert series of a 
weighted polynomial ring $k[x,y,z]$. 
\end{theorems} 
It is known that three dimensional Artin-Schelter regular algebras 
have all expected nice homological properties. For example they are 
both left and right noetherian domains. 

In this paper we restrict ourselves to Koszul three dimensional 
Artin-Schelter regular algebra. These have three generators and three 
defining relations in degree two. 
The 
minimal resolution of $k$ has the form 
\[ 
0 \r A(-3) \r A(-2)^{3} \r A(-1)^{3} \r A \r k_{A} \r 0 
\] 
hence $q_{k}(t) = (1-t)^{3}$ and the Hilbert series of $A$ is the same 
as that of the commutative polynomial algebra $k[x,y,z]$ with standard grading. 

Such algebras are also referred to 
as {\it quantum polynomial ring in three variables}. 
The corresponding $\Proj A$ will be 
called a {\it quantum projective plane} and will be denoted by 
$\P2q$. 

So let $A$ be a quantum polynomial ring in three variables. 
A \emph{linear module of dimension $d$} over $A$ is a cyclic 
$A$-module 
generated in 
degree zero with Hilbert series $(1-t)^{-d}$. Linear modules of dimension 
one and two are respectively called point and line modules.  The 
images of these objects in $\coh 
(\P2q)$  will be called point and line objects respectively. Line and 
point modules are classified in \cite{ATV1,ATV2}. 

Line modules are of the form $A/uA=L$ with $u\in A_1$. Hence 
line modules correspond naturally to lines in $\PP^2$.  To classify 
point modules we 
write the relations of $A$ as $f_{i} = \sum_{j=1}^{3} m_{ij}x_{j}$. 
Set $M = (m_{ij})_{i,j}$. We introduce auxiliary (commuting) 
variables $x_i^{(p)}$ (for $p\in \ZZ$)  and for a monomial 
$u=x_{i_0}\cdots x_{i_n}$ 
we define the \emph{multilinearization} of $m$ as $\widetilde{m}$ as 
$x_{i_0}^{(0)}\cdots x_{i_n}^{(n)}$. We extend this operation linearly to 
homogeneous polynomials in the variables $(x_i)_i$. 

Let 
$\Gamma \subset \PP^{2} \times \PP^{2}$ 
denote the locus of common zeros of the $\tilde{f_{i}}$. 
It turns out that $\Gamma$ is the graph of an automorphism $\sigma$ 
of $E=\pr_1(\Gamma)$, the locus of zeros of the multihomogenized polynomial 
$\det(\widetilde{M})$. If 
$\det(\widetilde{M})$ is not identically zero then $E$ is a 
divisor of degree $3$ in $\PP^{2}$. We then say that $A$ is {\it elliptic}. 
Otherwise, $E$ is all of $\PP^{2}$ and we call $A$ {\it linear} in this 
case. 

The connection between $E$ and point modules is as follows: let 
$P=\sum ke_u$ be a point module where $e_u\in P_u$. Put $e_u 
x_i=e_{u+1} \lambda_i^{(u)}$ with $ \lambda_i^{(u)}\in k$. From the fact 
that $e_0 f_i=0$ we deduce 
that $((\lambda_i^{(0)})_i,(\lambda_i^{(1)})_i)\in \Gamma$ and hence 
$(\lambda_i^{(0)})_i\in E$. This construction is reversible and defines 
a bijection between the closed points of $E$ and the point modules 
over $A$. If $P_q$ is the pointmodule corresponding to $q\in E$ then 
we have $P_q(1)_{\ge 0}=P_{\sigma p}$.

Let $j:E\r \PP^2$ be the inclusion and put $\Lscr=j^\ast \Oscr_{\PP^2}(1)$. 
Associated to the geometric data $(E,\sigma,\Lscr)$ is a so-called 
``twisted'' homogeneous coordinate ring $B=B(E,\sigma,\Lscr)$. This is 
a special case of a 
general construction in \cite{AZ}. See also \cite{AVdB}. Denote the 
auto-equivalence $\sigma_\ast(-\otimes_E \Lscr)$ by $-\otimes 
\Lscr_\sigma$. For $\Mscr\in \Qcoh(X)$ put 
$\Gamma_\ast(\Mscr)=\oplus_u \Gamma(E,\Mscr\otimes 
(\Lscr_\sigma)^{\otimes u})$ and $B=\Gamma_\ast(\Oscr_E)$. It is easy 
to see that $B$ has a natural ring structure and $\Gamma_\ast(\Mscr)$ 
is a right $B$-module. A straightforward verification shows 
\[ 
\Gamma_\ast(\Oscr_q)=P_q. 
\] 
In \cite{ATV2} it is shown that there is a surjective morphism 
$p:A\r B$ of graded $k$-algebras.  Its kernel is trivial in 
the linear case and it is generated by a regular normalizing element 
$g$ of degree three in the elliptic case.  All point modules are 
$B$-modules. In other words: $g$ annihilates all point modules. 

By analogy with the commutative case we may say that $\Proj A$ contains $\Proj 
B$ as a ``closed'' subscheme. Though the structure of $\Proj A$ is 
somewhat obscure, that of $\Proj B$ is well understood. 

Indeed it follows from \cite{AZ,AVdB} that 
the functor $\Gamma_\ast:\Qcoh(E)\r \GrMod(B)$ defines an equivalence 
$\Qcoh(E)\cong \Tails(B)$. The inverse of this equivalence and its 
composition with $\pi:\GrMod(B)\r \Tails(B)$ are both denoted by 
$\widetilde{(-)}$ 

For further properties of point modules and line modules over three 
dimensional quantum polynomial algebras we refer to \cite{ATV1,ATV2}. 

We will frequently use the following result 
\begin{lemmas} \label{ref-2.8.2-18} Assume that we are in the elliptic 
case. Let $M\in 
\grmod(A)$ be such that $M/Mg\in \tors 
A$. Then $\gkdim M=1$. If $\sigma$ has infinite order then $M\in 
\tors(A)$. 
\end{lemmas} 
\begin{proof} Multiplication induces an isomorphism $M_n\cong M_{n+3}$ 
for large $n$. Hence $\gkdim M=1$. Furthermore $(M_g)_0$ is a finite 
dimensional representation of $(A_g)_0$. It is is shown in 
\cite{ATV2} that if $\sigma$ has infinite order then $(A_g)_0$ is a 
simple ring. In particular it has no finite dimensional 
representations. Thus $(M_g)_0=0$. This implies $M\in \tors(A)$. 
\end{proof}

In the sequel it will be useful to cast the relationship between the 
noncommutative graded ring $A$ and the commutative scheme $E$ into the 
language of 
noncommutative algebraic geometry exhibited in \cite{smith1,VdB19} 
although we will 
use this language only in an intuitive way. Let $X=\Proj A$, $Y=\Proj 
B$. 

We define a map of noncommutative schemes $i:E\r X$ by 
\begin{align*} 
i^\ast \pi M&=(M\otimes_A B)\tilde{}\\ 
i_\ast\Mscr&=\pi(\Gamma_\ast(\Mscr)_A) 
\end{align*} 

We will call $i^\ast(\pi M)$ the \emph{restriction} of $\pi M$ to 
$E$. $i_\ast$ is clearly an exact functor. For the left derived 
functor of $i^\ast$ we have: 
\begin{lemmas} \label{ref-2.8.3-19} If $M\in D^-(\GrMod(A))$ then 
$Li^\ast(\pi M)=(M\LotimesA B)\tilde{}$ 
\end{lemmas} 
\begin{proof} One shows first that the objects $\pi F$ where $F$ is a 
finitely generated graded free $A$-module are acyclic for $i^\ast$ in the 
sense of \cite{RD}. Then the lemma follows by replacing $M$ 
by a resolution of finitely generated free $A$-modules. 
\end{proof} 
We easily obtain the following consequence: 
\begin{lemmas} \label{ref-2.8.4-20} Assume that we are in the elliptic case and 
let $\Mscr\in D^-(\Qcoh(\P2q))$. Then there are short exact sequences: 
\[ 
0\r i^\ast H^j(\Mscr)\r H^j(Li^\ast \Mscr)\r L_1i^\ast H^{j+1}(\Mscr)\r 0
\] 
\end{lemmas} 
\begin{proof} Take $M\in D^-(\GrMod(A))$ such that $\Mscr=\pi M$. We may 
assume that $M$ is given by a right bounded complex of graded 
projective $A$-modules. The lemma now follows 
by applying $\pi$ to the long exact homology 
sequence associate to the short exact sequence of complexes 
\[ 
0\r Mg\r M \r M/Mg\r 0\qed 
\] 
\def\qed{}\end{proof} 
\subsection{Three dimensional Sklyanin algebras} \label{ref-2.9-21} 
Below we are interested in Sklyanin algebras of dimension 
three which are elliptic Artin-Schelter regular algebras such that the 
corresponding elliptic curve $E$ is smooth and the automorphism is a 
translation. More specificly, we are 
interested in the algebras 
\[ 
\Skly_{3}(a,b,c) = k\{ x,y,z \} / (f_{1}, f_{2}, f_{3}) 
\] 
where $f_{1},f_{2},f_{3}$ are the quadratic equations 
\begin{equation} \label{ref-2.4-22} 
\left\{ 
\begin{array}{l} 
f_{1} = ayz + bzy + cx^{2}  \\ 
f_{2} = azx + bxz + cy^{2}\\ 
f_{3} = axy + byx + cz^{2} 
\end{array} \right. 
\end{equation} 
and $(a,b,c) \in \PP^{2}\setminus F$ where 
\[ 
F = \{ (a,b,c) \in \PP^{2} \mid abc = 0 \mbox{ or } 
a^{3} = b^{3} = c^{3} \mbox{ or } (3abc)^{3} = 
(a^{3} + b^{3} + c^{3})^{3} \}. 
\] 
The algebras $\Skly_{3}(a,b,c)$ are elliptic quantum polynomial 
rings. They correspond to Artin-Schelter algebras 
of dimension three where, in the associated geometric data, $E$ is a smooth 
elliptic curve and $\sigma$ 
is given by translation under the group law. We refer to \cite{ATV1} 
for the description of $E$ and $\sigma$. The regular normalizing 
element $g$ of degree three turns 
out to be central in this case. 

Put $A=\Skly_{3}(a,b,c)$. 
Combining the results in \cite{VdB16} 
with Theorem \ref{ref-2.3.1-8} we see that Serre duality for $A$ 
takes a particularly simple form: 
\begin{theorems} \label{ref-2.9.1-23} Let $\Mscr,\Nscr\in 
D^b(\tails(A))$. Then there are natural isomorphisms 
\[ 
\Ext^{i}(\Mscr,\Nscr)\cong \Ext^{2-i}(\Nscr,\Mscr(-3))^\ast 
\] 
\end{theorems} 
\begin{corollarys} Let $\Mscr\in D^b(\tails(A))$ and let $\Pscr\in \tails(A)$ 
be a point object corresponding to $p\in E$. Then 
\[ 
\Ext^i(\Pscr,\Mscr)\cong \Ext^{2-i}(\Mscr,\Pscr')^\ast  
\] 
where $\Pscr'$ is the point object corresponding to $\sigma^{-3} 
p$. 
\end{corollarys} 
\section{Cohomology of rank one sheaves on a quantum projective plane} 
In this section, $A$ will be a quantum polynomial ring in three 
variables, and $\P2q = \Proj A$ the associated quantum projective 
plane. As usual $\Oscr=\pi A$.     

We say that a graded right $A$-module $M \neq 0$ is {\it torsion} if $\rk M = 
0$. $M$ is called {\it torsion-free} if $M$ contains no torsion 
submodule. This is the same as saying that $M$ is pure 
three dimensional. We use the same terminology for objects in $\coh(\P2q)$. 

The graded right ideals of $A$ are, up to isomorphism, precisely the 
shifts of 
torsion-free rank one 
right $A$-modules. 

A torsion-free rank one graded 
$A$-module $I$ gives rise to a torsion-free coherent sheaf 
$\Iscr = \pi I$ on $\P2q$ of rank one. 
Conversely, every torsion-free $\Iscr \in \coh(\P2q)$ 
determines a torsion-free rank one graded 
$A$-module $\omega \Iscr$. 

Any shift $l$ of a torsion-free rank one graded 
$A$-module $I$  gives rise 
to a torsion-free rank one coherent sheaf $\Iscr(l) = \pi I(l)$ on $\P2q$. 
Our  first aim is to normalize this shift. \\ 

We will use the following natural basis for $\K0$. 
\begin{proposition} \label{ref-3.1-24} 
Let $P$ be a point module and $S$ a line module over $A$. Denote the 
corresponding objects  in $\coh(\P2q)$ by $\Pscr$ and $\Sscr$. \\ 
Then 
$\{ [\Oscr],[\Sscr],[\Pscr]\}$ is a $\ZZ$-module basis of $\K0$, which 
does not depend on the particular choice of $S$ and $P$, and 
the action of the shift functor on that basis is 
\begin{equation} \label{ref-3.1-25} 
\begin{aligned} 
\left[\Oscr(1) \right] & =  [\Oscr] + [\Sscr] + [\Pscr] \\ 
\left[\Sscr(1) \right] & =  [\Sscr] + [\Pscr] \\ 
\left[\Pscr(1) \right] & =  [\Pscr] 
\end{aligned} 
\end{equation} 
\end{proposition} 
\begin{proof} 
It follows from Theorem \ref{ref-2.5.1-11} that the class in $K_0(\P2q)$ of 
an object $\pi M$ depends only on the Hilbert series of $M$. Thus 
$[\Sscr]$ and $[\Pscr]$ are indeed independent of  the particular 
choice of $S$ and $P$. 

Using a computation with Hilbert series we see that the images of 
$[\Oscr],[\Sscr]$ and $[\Pscr]$ 
under the isomorphism $\theta$ of Theorem \ref{ref-2.5.1-11} 
\[ 
\theta : K_{0}(\P2q) \r \ZZ[t,t^{-1}]/ (1-t)^{3} 
\] 
are respectively $\overline{1}$, $\overline{1-t}$, 
$\overline{(1-t)^2}$. Furthermore the shift functor 
corresponds to multiplication by $t^{-1}$.  This easily yields what we want. 
\end{proof} 
From now on, we fix such a $\ZZ$-module basis $\{ [\Oscr],[\Sscr],[\Pscr]\}$ 
of 
$\K0$. 
For any coherent sheaf $\Jscr$ on 
$\P2q$ we may write 
\[ 
[\Jscr] = r[\Oscr] + a[\Sscr] + b[\Pscr] 
\] 
where $r$ is the rank of $\Jscr$. \\ 
It follows from (\ref{ref-3.1-25}) that we have 
\begin{eqnarray} \label{ref-3.2-26} 
[\Jscr(l)] = r[\Oscr] + (a + lr)[\Sscr] + (\frac{1}{2}l(l+1)r + la + b)[\Pscr] 
\end{eqnarray} 
for all integers $l$. 
\begin{proposition} \label{ref-3.2-27} 
\begin{enumerate} 
\item 
Let $\Iscr$ be a coherent sheaf on $\P2q$ of rank one, and write 
$[\Iscr] = [\Oscr] + a[\Sscr] + b[\Pscr]$. 
Then there is an unique shift $c$ (namely $-a$) and an integer $n$ such that 
\[ 
[\Iscr(c)] = [\Oscr] - n[\Pscr]. 
\] 
Moreover, $n = \frac{1}{2}a(a + 1) - b$. 
\item 
Let $\Fscr$ be a coherent sheaf on $\P2q$ of rank zero, and write 
$[\Fscr] = u[\Sscr] + v[\Pscr]$. Then $u=e_1(\Fscr)$. If $u = 0$, then 
$v=e_0(\Fscr)$. 
\end{enumerate} 
\end{proposition} 
\begin{proof} 
For the first part, use (\ref{ref-3.2-26}). The uniqueness is easy to see.\\ 
For the second statement, take the image of 
$[\Fscr] = u[\Sscr] + v[\Pscr]$ under the 
isomorphism 
$\theta$ of Theorem \ref{ref-2.5.1-11}. Take $F$ such that  $\Fscr=\pi F$. 
We obtain 
\[ 
q_{F}(t) = uq_{L}(t) + vq_{P}(t) + f(t)q_{k}(t) 
\] 
for a suitable $f(t) \in \ZZ[t,t^{-1}]$. Multiplying both sides 
with $h_{A}(t) = q_{k}(t)^{-1}$ yields (see (\ref{ref-2.3-10})) 
\[ 
h_{F}(t) = uh_{L}(t) + vh_{P}(t) + f(t) 
\] 
We find $e_2(F)=\lim_{t\r 1} (1-t)^2h_F(t)=u$ and if $u=0$ then 
$e_1(F) = \linebreak
\lim_{t\r 1} (1-t)h_F(t)=v$. 
\end{proof} 
We call the integer $n$ appearing in {Proposition} \ref{ref-3.2-27} the 
``invariant'' of $\Iscr$ (or of the corresponding  torsion-free rank 
one graded 
$A$-module $I$  such that 
$\Iscr=\pi I$). 
Note that two torsion-free rank one graded 
$A$-modules $I, J$ have the same 
invariant if and only if $\dim_{k} I_{i} = \dim_{k} J(d)_{i}$ for $i \gg 0$
and for a fixed integer $d$. 

We will call a torsion-free rank one coherent sheaf $\Iscr$ on $\P2q$ 
{\it normalized} if $[\Iscr] = [\Oscr] - n[\Pscr]$ for an integer $n$. 
We will prove later that this $n$ is actually positive. 
 
We will call a torsion-free reflexive rank one sheaf on $\P2q$ a 
\emph{line bundle}. Our aim is to classify line bundles on $\P2q$ up 
to shift. By the above discussion this is equivalent to classifying 
normalized line bundles up to isomorphism. 

It is also easy to see that through the functors $\pi$ and 
$\omega$ classifying line bundles up to shift is equivalent to 
classifying reflexive torsion-free rank one graded 
$A$-modules, also up to shift. 

We recall two elementary lemmas. 
\begin{lemma} \label{ref-3.3-28} 
Let $\Mscr, \Nscr$ be torsion-free coherent sheaves on $\P2q$ of rank one. 
Then every nonzero morphism in $\Hom(\Mscr,\Nscr)$ is injective. 
\end{lemma}
\begin{proof} 
$\Mscr$ and $\Nscr$ are critical of the same dimension. It is 
well-known that this implies that any map between them must be 
injective \cite{ATV2}. 
\end{proof}
\begin{lemma} 
\label{ref-3.4-29} 
Let $\Mscr \in \coh(\P2q)$. 
Then $\Mscr$ is reflexive
if and only if $\Mscr$ is torsion-free and $\Ext^1(\Nscr,\Mscr)=0$
for all $\Nscr \in \coh(\P2q)$ of dimension zero. 
\end{lemma}
\begin{proof}
Assume that $\Mscr$ is a reflexive coherent sheaf on $\P2q$.
By \eqref{ref-2.1-6} we need to prove the corresponding 
statement for $\grmod(A)$. Thus assume that $M$ is a reflexive $A$-module and 
$\gkdim N\le 1$. Assume that there is a non-split exact sequence 
\begin{equation} 
\label{ref-3.3-30} 
0\r M\r M'\r N\r 0 
\end{equation} 
By \cite[Theorem 4.1]{ATV2} one has $\underline{\Ext}^1(N,A)=0$. Hence 
we obtain $M^{\prime \ast}=M^\ast$ and thus 
$M=M^{\ast\ast}=M^{\prime\ast\ast}$. Thus the composition of $M\r M'\r 
M^{\prime\ast\ast}$ is an isomorphism, implying that the first map 
splits. This contradicts the non-triviality of the extention
\eqref{ref-3.3-30}. \\
For the other implication, let $\Mscr \in \coh(\P2q)$ be torsion-free and 
$\Ext^1(\Nscr,\Mscr)=0$ for all $\Nscr \in \coh(\P2q)$ of dimension zero. 
$M = \omega \Mscr$ is pure and $\gkdim M = 3$ since $\Mscr$ is pure two dimensional. 
By \cite[Corollary 4.2]{ATV2}
there is a canonical map $\mu : M \r M^{**}$ and $\ker \mu$ is the 
maximal submodule of $M$ which has $\gkdim < 3$. 
Hence $\mu$ is injective, and we have an exact sequence
\begin{equation} \label{mu}
0 \r M \r M^{**} \r \coker \mu \r 0
\end{equation}
where $\gkdim(\coker \mu) \leq 1$.
Applying $\pi$ on (\ref{mu}) yields 
\begin{equation} \label{pimu}
0 \r \Mscr \r \pi M^{**} \r \Nscr \r 0
\end{equation}
where $\Nscr = \pi \coker \mu$.
Now $\Nscr$ must be zero, otherwise $\dim \Nscr = 0$ and since
$\Ext^{1}(\Nscr, \Mscr) = 0$ the sequence (\ref{pimu}) would split,
which is impossible because $M^{**}$ is pure three dimensional. 
Hence $\Mscr = \pi M = \pi M^{**}$ and thus $\Mscr$ is reflexive.
\end{proof}
Now we can partially compute the cohomology of line bundles on $\P2q$. 
This computation is similar to the one for the homogenized Weyl 
algebra in \cite{lebruyn1}. However the computations for the 
homogenized Weyl algebra rely  on the existense of a central 
element in degree one. So they do not apply in a 
straightforward way to the case we 
consider. 
\begin{theorem}  \label{ref-3.5-31} 
Let $\Mscr$ be a rank one torsion-free coherent sheaf  on $\P2q$ 
where \linebreak
$[\Mscr] = [\Oscr] - n[\Pscr]$. 
Assume that 
$\Mscr \ncong \Oscr$. 
Then 
\begin{enumerate} 
\item 
$H^{0}(\P2q , \Mscr (l)) = 0 \mbox{ for } l \leq 0, $\\ 
$H^{2}(\P2q , \Mscr (l)) = 0 \mbox{ for } l \geq -2;$ 
\item 
$\chi (\Oscr,\Mscr(l)) = \frac{1}{2}(l+1)(l+2) - n \mbox{ for all } 
l \in \ZZ;$ 
\item 
$\dim_{k} H^{1}(\P2q,\Mscr) = n - 1$ \\ 
$\dim_{k} H^{1}(\P2q,\Mscr(-1)) = n$ \\ 
$\dim_{k} H^{1}(\P2q,\Mscr(-2)) = n$ 
\item 
$H^j(\P2q,\Mscr)=0$ for $j\ge 3$. 
\end{enumerate} 
As a consequence, $n$ is positive and nonzero. 

If $\Mscr$ is a line bundle then we have in addition: 
$H^{2}(\P2q , \Mscr (l)) = 0 \mbox{ for } l=-3$ and $\dim_{k} 
H^{1}(\P2q,\Mscr(-3)) = n - 1$. 
\end{theorem} 
\begin{proof} 
That $H^j(\P2q,\Mscr)=0$ for $j\ge 3$ is part of Theorem \ref{ref-2.6.2-13}. 

To prove the rest of the current theorem we 
first let $l \leq 0$. Suppose $f$ is a nonzero morphism in 
$\Hom (\Oscr, \Mscr(l))$. By   lemma \ref{ref-3.3-28} 
$f$ is injective and from the exact sequence 
\begin{eqnarray} \label{ref-3.4-32} 
0 \r \Oscr \r \Mscr(l) \r \coker f \r 0 
\end{eqnarray} 
we get $[\coker f] = l[\Sscr] + (l(l+1)/2 - n)[\Pscr]$. 
Using Proposition \ref{ref-3.2-27} gives $l \geq 0$, thus $l = 0$ and 
$[\coker f] = -n[\Pscr]$. Hence by the discussion in 
\S\ref{ref-2.7-15} together with Proposition 
\ref{ref-3.2-27} we obtain $\dim \coker f=0$. 
 By 
lemma \ref{ref-3.4-29} $\Ext^{1}(\coker f, \Oscr) = 0$. This means that 
the exact sequence 
(\ref{ref-3.4-32}) splits hence $\Mscr$ is not torsion-free. A 
contradiction. We conclude that $\Hom(\Oscr,\Mscr(l)) = 0$ for $l \leq 
0$. \\ 
Second, let $l \geq -2$. Serre duality (Theorem \ref{ref-2.9.1-23}) 
yields 
\[ 
\Ext^{2}(\Oscr,\Mscr(l))^{*} \cong \Hom(\Mscr(l+3),\Oscr). 
\] 
If $g$ is a nonzero morphism in $\Hom(\Mscr(l+3),\Oscr)$ then $g$ is 
injective, and from the exact sequence 
\begin{equation} 
\label{ref-3.5-33} 
0 \r \Mscr(l+3) \r \Oscr \r \coker g \r 0 
\end{equation} 
we get $[\coker g] = u[\Sscr] + v[\Pscr]$ where $u = -(l+3)$ and 
$v = n - (l+3)(l+4)/2$. By Proposition \ref{ref-3.2-27} $u\ge 0$ but $l\ge 
-2$ implies 
$u<0$. This yields a contradiction. 

Assume now  $l \ge -3$ and $\Mscr$ reflexive. By the same reasoning as 
above we 
obtain $l=-3$ and thus the dimension of 
$\omega\coker g$ is zero. By lemma \ref{ref-3.4-29} it follows that 
\eqref{ref-3.5-33} splits. But this contradicts the fact that $\Oscr$ is 
torsion-free. 

For the second part we use Theorem \ref{ref-2.6.2-13} to obtain 
\[ 
\chi(\Oscr,\Oscr(l)) = \frac{1}{2}(l+1)(l+2) \mbox{ for all } l \in \ZZ 
\] 
and from Proposition \ref{ref-3.1-24} we deduce 
\begin{eqnarray*} 
[\Sscr] & = & -[\Oscr(2)] + 3[\Oscr(1)]  - 2[\Oscr] \\ 
\left[\Pscr \right] & = & \hspace{0.25cm}
[\Oscr(2)] - 2[\Oscr(1)]  + \hspace{0.25cm}
[\Oscr] 
\end{eqnarray*} 
Combining these results yields $\chi(\Oscr,\Sscr) = 1$ and 
$\chi(\Oscr,\Pscr) = 1$. 
Now we use (\ref{ref-3.2-26}) to obtain 
\begin{equation*} 
\begin{split} 
\chi(\Oscr,\Mscr(l)) & = \chi(\Oscr, \Oscr) +l \chi(\Oscr,\Sscr)  + 
 \left(\frac{1}{2}l(l+1) - n\right) \chi(\Oscr,\Pscr)\\ 
& = \frac{1}{2}(l+1)(l+2) - n 
\end{split} 
\end{equation*} 

Finally, we combine the first two results of the theorem. 
If $-2 \leq l \leq 0$ (or $-3\leq l \leq 0$ if $\Mscr$ is reflexive) the 
first statement gives 
\begin{equation*} 
\begin{split} 
\chi(\Oscr,\Mscr(l)) & = \dim_{k}H^{0}(\P2q,\Mscr(l)) - 
\dim_{k}H^{1}(\P2q,\Mscr(l)) 
+ \dim_{k}H^{2}(\P2q,\Mscr(l)) \\ 
& = - \dim_{k}H^{1}(\P2q,\Mscr(l)) 
\end{split} 
\end{equation*} 
and comparing with the expression 
$\chi (\Oscr,\Mscr(l)) = \frac{1}{2}(l+1)(l+2) - n$ 
completes the proof. 
\end{proof} 
Using Theorem \ref{ref-3.5-31} the 
torsion-free rank one graded 
$A$-modules having invariant zero 
are easy to determine. 
\begin{corollary} 
\label{ref-3.6-34} 
Let $\Iscr$ be a torsion-free coherent sheaf of rank one on $\P2q$ 
with invariant $n$. 
Then 
\[ 
n = 0 \Leftrightarrow \Iscr\cong \Oscr(d) 
\] 
for some integer $d$. 
\end{corollary} 
\begin{proof} 
If $\Iscr\cong \Oscr(d)$ then clearly $n=0$. Assume conversely 
$n=0$. We may assume that $\Iscr$ is normalized. If $\Iscr\not\cong 
\Oscr$ then by Theorem \ref{ref-3.5-31} $n>0$. Since $n=0$ we obtain 
$\Iscr\cong\Oscr$ by contraposition. 
\end{proof} 

\section{Restriction of coherent sheaves} 
In this section, $A$ will be a Sklyanin algebra $\Skly_{3}(a,b,c)$ as 
defined in \S\ref{ref-2.9-21}. We recycle the notations of sections 
\S\ref{ref-2.6-12}-\S\ref{ref-2.9-21}.  In particular the symbols 
$\Oscr,E,\sigma,\Lscr,B,i$ have their usual meaning. 

Note that $E$ is 
a smooth elliptic curve. We fix a grouplaw on $E$. Then $\sigma$ is a 
translation by some element $\xi \in E$. 

The dimension of 
objects in $\grmod(B)$ or $\tails(B)$ will be computed in $\grmod(A)$ or 
$\tails(A)$. The dimension of objects in $\coh(E)$ is the dimension of 
their support. 

There is a group homomorphism 
\[ 
K_0(\P2q)\r K_0(E):[\Mscr]\mapsto [i^\ast\Mscr]-[L_1i^\ast \Mscr] 
\] 
which as usual is also denoted by $i^\ast$. 
\begin{lemma} 
\label{ref-4.1-35} 
We have 
\begin{align*} 
i^\ast[\Oscr]&=[\Oscr_E]\\ 
i^\ast[\Sscr]&=[\Oscr_u]+ [\Oscr_v]+ [\Oscr_w]&&\text{$u,v,w$ arbitrary but 
colinear}\\ 
i^\ast[\Pscr]&=[\Oscr_p]-[\Oscr_{p^{\sigma^{-3}}}] &&\text{$p$ arbitrary} 
\end{align*} 
\end{lemma} 
\begin{proof} This follows easily from lemma \ref{ref-2.8.3-19} 
\end{proof} 
According to \cite[Ex II. 6.11]{H} we have $K_0(E)\cong \ZZ\oplus 
\Pic(E)$. The projection $K_0(E) \r \ZZ$ is given by the rank and the 
projection $K_0(E) \r \Pic(E)$ is given the first Chern class. If 
$\Escr$ is a vector bundle on $E$ then $c_1(\Escr)=\wedge^{\rk 
\Escr}\Escr$. We also have for $q\in E$: $c_1(\Oscr_q)=\Oscr_E(q)$.

There is a homomorphism $\deg:\Pic(E)\r \ZZ$ which assigns to a line 
bundle its degree. For simplicity we will denote the composition 
$\deg\circ\, c_1$ also by $\deg$. If $\Uscr$ is a line bundle then $\deg 
[\Uscr]=\deg\Uscr$. If $F\in\coh(E)$ has finite length then $\deg 
[F]=\length F$ \cite[Ex. 6.12]{H}. From lemma \ref{ref-4.1-35} we deduce 
that if $[\Mscr]=a[\Oscr]+b[\Sscr]+c[\Pscr]$ then 
\begin{gather} 
\rk i^\ast [\Mscr]=a=\rk \Mscr \label{ref-4.1-36}\\ 
\deg i^\ast [\Mscr]=3b\label{ref-4.2-37} 
\end{gather} 

\begin{lemma} \label{ref-4.2-38} 
\begin{enumerate} 
\item 
If $M\in \grmod(B)$ is pure two dimensional then 
$\tilde{M}\in \coh(E)$ is pure one dimensional. 
\item If $\Nscr \in \coh(E)$ is pure one dimensional then 
$\Gamma_\ast(\Nscr)$ is pure two dimensional. 
\end{enumerate} 
\end{lemma} 
\begin{proof}  The indecomposable objects in $\coh(E)$ are vector 
bundles and finite length objects. 
Using Riemann-Roch it is easy to see that if $0\neq 
\Uscr\in \coh(E)$ then $\gkdim \Gamma_\ast(\Uscr)=\dim 
\Uscr+1$. From this we deduce that if $V\in \grmod(B)$ is not in 
$\tors(B)$ then $\gkdim V=\dim \tilde{V}+1$. The lemma now easily 
follows. 
\end{proof} 
We deduce 
\begin{proposition}  \label{ref-4.3-39} 
\begin{enumerate} 
\item 
If $\Mscr\in \coh(\P2q)$ is reflexive then $i^\ast 
\Mscr$ is a vector bundle on $E$ and $L_{j}i^\ast \Mscr=0$ for $j>0$. 
\item If $\Mscr$ is a line bundle then so is $i^\ast \Mscr$. 
\item If $\Mscr$ is a line bundle then $\Mscr$ is normalized if and 
only if $\deg i^\ast\Mscr=0$. 
\item If $\Mscr$ is a normalized line bundle with invariant $n$ then 
\[ 
c_1(i^\ast\Mscr)=\Oscr((o)-(3n\xi)) 
\] 
where ``$o$'' is the origin for the group law. 
\end{enumerate} 
\end{proposition} 
\begin{proof} 
\begin{enumerate} 
\item 
We have $\Mscr=\pi M$ where $M$ is reflexive. In 
particular $M$ is torsion-free.  By lemma \ref{ref-2.8.3-19} it follows 
that $L_{j}i^\ast \Mscr=0$ for $j>0$ and \linebreak
$i^\ast \Mscr=(M/Mg)\tilde{}$. 

If $M/Mg$ contains a nonzero submodule 
$N/Mg$ of GK-dimension $\le 1$ then $N$ represents an element of 
$\Ext^1(N/Mg, Mg)$ which must be zero by lemma \ref{ref-3.4-29} (or 
rather its proof). Thus $N/Mg\subset N\subset M$. This is impossible 
since $M$ is torsion-free. 

Hence $M/Mg$ is pure of GK-dimension 2.  By the previous lemma it 
follows that $(M/Mg)\tilde{}$ is a vector bundle. 
\item 
This follows from \eqref{ref-4.1-36}. 
\item This follows from \eqref{ref-4.2-37}. 
\item We have $[\Mscr]=[\Oscr]-n[\Pscr]$. By lemma \ref{ref-4.1-35} we 
obtain 
$[i^\ast\Mscr]=[\Oscr_E]-n[\Oscr_p]+n[\Oscr_{p^{\sigma^{-3}}}]$. 
Hence $c_1 (i^\ast \Mscr)=\Oscr(n(p^{\sigma^{-3}})-n(p))$. Now 
$n(p^{\sigma^{-3}})-n(p)$ and $(o)-(3n\xi)$ are both divisors of 
degree zero which have the same sum for the group law. Hence they 
are linearly equivalent by \cite[IV Thm 4.13B ]{H}. This finishes 
the proof. \qed 
\end{enumerate} 
\def\qed{}\end{proof} 
Now we prove a converse of Proposition  \ref{ref-4.3-39}. 
\begin{proposition}\label{ref-4.4-40} 
Assume that $\sigma$ has infinite order and that 
$\Mscr\in D^b(\coh(\P2q))$ is such that $Li^\ast \Mscr$ is a vector 
bundle on $E$. Then $\Mscr$ is a reflexive object in $\coh(\P2q)$. 
\end{proposition} 
\begin{proof} 
It follows from lemma \ref{ref-2.8.4-20} that $i^\ast H^j(\Mscr)=0$ for $j\neq 
0$. Then it follows from lemma \ref{ref-2.8.2-18} that 
$\Mscr\in\coh(\P2q)$ and $L_1 i^\ast \Mscr=0$, using lemma 
\ref{ref-2.8.4-20} again. 

Pick an object $M$ in $\grmod(A)$ such that $\pi 
M=\Mscr$. We may assume that $M$ contains no subobject in $\tors(A)$.  By 
lemma \ref{ref-2.8.3-19} we have $ L_1 i^\ast 
\Mscr=\ker(M(-3)\xrightarrow{\times g} M)\tilde{}$. Thus 
$\ker(M(-3)\xrightarrow{\times g} M) \in \tors(A)$. Since $M$ contains 
no subobject in $\tors(A)$ it follows that $M$ is $g$-torsion 
free. Furthermore by lemma \ref{ref-4.2-38} 
$\Gamma_\ast(i^\ast\Mscr)=\Gamma_\ast((M/Mg)\tilde{})$ is pure 
two dimensional. If $T$ is the maximal submodule of $M/Mg$ which is in 
$\tails(A)$ then since $(M/Mg)/T\subset \Gamma_\ast((M/Mg)\tilde{})$ we 
obtain that \linebreak
$(M/Mg)/T$ is pure two dimensional. 

We now claim that $M$ is pure three dimensional. Let $N$ be the 
maximal submodule of $M$ of dimension $\le 2$. Then $C=M/N$ is pure 
three dimensional and in particular $g$-torsion-free. Hence we have a 
short exact sequence 
\[ 
0\r N/Ng\r M/Mg \r C/Cg\r 0 
\] 
By the purity of $(M/Mg)/T$ it follows that $N/Ng\subset T$ and hence 
$N/Ng \in \tails(A)$. It follows from lemma \ref{ref-2.8.2-18} than $N\in 
\tails(A)$ and hence $N=0$. This shows that $M$ is pure. 

Put $Q=M^{\ast\ast}/M$.  Thus we obtain an exact sequence 
\[ 
0\r \Tor^A_1(Q,B)\tilde{}\r (M\otimes_A B)\tilde{}\r 
 (M^{\ast\ast}\otimes_A B)\tilde{}\r (Q\otimes_A B)\tilde{}\r 0 
\] 
By \cite{ATV2} we have $\gkdim Q\le 
1$. Thus we have $\gkdim \Tor^A_1(Q,B)\le \gkdim Q\le 1$. So by the proof 
of lemma 
\ref{ref-4.2-38}, $\dim \Tor^A_1(Q,B)\tilde{}\le 
0$. Since $(M\otimes_A B)\tilde{}$ is a vector bundle by hypotheses it 
contains no finite dimensional subobjects and we 
obtain $\Tor^A_1(Q,B)\tilde{}=0$. Thus $\Tor^A_1(Q,B)\in \tors(A)$. 
Thus, in high degree, multiplication by $g$ is an isomorphism on 
$Q$. But then by lemma \ref{ref-2.8.2-18} $Q\in \tors(A)$. Hence 
$\Mscr=\pi M=\pi M^{\ast\ast}$ and thus $\Mscr$ is reflexive. 
\end{proof}

\section{Elliptic quantumspaces} 
\subsection{Generalities} 
Let $A$ be a Sklyanin algebra $\Skly_{3}(a,b,c)$. We use again our 
standard notations as in the previous section. 

We set $\Escr = \Oscr(2) \oplus \Oscr(1) \oplus \Oscr$ and ${{D}} = 
\Hom_{\P2q}(\Escr, \Escr) = \bigoplus_{i,j = 
0}^{2}{\Hom_{\P2q}(\Oscr(i),\Oscr(j))}$ the algebra of endomorphisms 
of $\Escr$. We consider the left exact functor 
$\Hom_{\P2q}(\Escr,-)$ which takes coherent sheaves on $\P2q$ to 
right ${{D}}$-modules. 

$\Hom_{\P2q}(\Escr,-)$ extends to a functor $\RHom_{\P2q}(\Escr,-)$ 
on bounded derived categories 
\begin{equation} 
\label{ref-5.1-41} 
\RHom_{\P2q}(\Escr,-) : \D^{b}(\coh(\P2q)) \r \D^{b}(\mod(D)). 
\end{equation} 
This is done as follows: $\Qcoh(\P2q)$ has enough 
injectives 
and this yields a functor  $\RHom_{\P2q}(\Escr,-):
\D^b_{\coh(\P2q)}(\Qcoh(\P2q))\r D^b_{\mod(D)}(\Mod(D))$. 
Now $\coh(\P2q)$ and $\mod(D)$ are noetherian 
abelian categories and this yields equivalences \linebreak 
$\D^{b}(\coh(\P2q))\cong \D^b_{\coh(\P2q)}(\Qcoh(\P2q))$ and 
$\D^{b}(\mod(D))\cong D^b_{\mod(D)}(\Mod(D))$ (lemma \ref{ref-2.2-5}). 
The functor 
\eqref{ref-5.1-41} is obtained by composing with these equivalences. 

In a similar way as in  \cite[Theorem 6.2]{Bondal1} one shows that 
$\RHom_{\P2q}(\Escr,-)$ 
is an equivalence of 
derived categories. 
The inverse functor is given by $- \LotimesD \Escr$. 
For a non-negative integer $i$ the equivalence restricts to an 
equivalence between $\Xscr_{i}$ and $\Yscr_{i}$ where 
$\Xscr_{i} \subset \coh(\P2q)$ is the full subcategory with objects 
\[ 
\Xscr_{i} = \{ \Mscr \in \coh(\P2q) \mid \Ext_{\P2q}^{j}(\Escr,\Mscr) = 0 
\mbox{ for } j \neq i \} 
\] 
and $\Yscr_{i} \subset \mod(D)$ the full subcategory with objects 
\[ 
\Yscr_{i} = \{ M \in \mod(D) \mid \Tor_{j}^{{{D}}}(M,\Escr) = 0 
\mbox{ for } 
j \neq i \}. 
\] 
The inverse equivalences between these categories are given by 
$\Ext_{\P2q}^{i}(\Escr, -)$ and $\Tor_{i}^{{{D}}}(-, \Escr)$. 

\smallskip 

Let $(\Delta,R)$ be the quiver 
\begin{eqnarray} \label{ref-5.2-42} 
\begin{array}{ccccc} 
\quad & {\stackrel{X_{-2}}{\lr}} & \quad & {\stackrel{X_{-1}}{\lr}} & 
\quad \\ 
-2 & {\stackrel{Y_{-2}}{\lr}} & -1 & {\stackrel{Y_{-1}}{\lr}} &  0 \\ 
\quad & {\stackrel{Z_{-2}}{\lr}} & \quad & {\stackrel{Z_{-1}}{\lr}} & 
\quad 
\end{array} 
\end{eqnarray} 
with relations 
\begin{eqnarray} \label{ref-5.3-43} 
\left\{ 
\begin{array}{l} 
aY_{-2}Z_{-1} + bZ_{-2}Y_{-1} + cX_{-2}X_{-1}  =  0 \\ 
aZ_{-2}X_{-1} + bX_{-2}Z_{-1} + cY_{-2}Y_{-1}  =  0 \\ 
aX_{-2}Y_{-1} + bY_{-2}X_{-1} + cZ_{-2}Z_{-1}  =  0 
\end{array} 
\right. 
\end{eqnarray}

We 
write $\Mod(\Delta)$  for the  category of representations of the 
quiver $\Delta$ (representations are always assumed to satisfy the 
relations \eqref{ref-5.3-43}). If $i=-2,-1,0$ then we denote by 
$P_i,S_i$ respectively the 
projective representation and the simple representation 
corresponding to $i$. 

It is easy to see that ${{D}}\cong 
k\Delta/(R)$. Since the category $\Mod(\Delta)$ of 
representations of $\Delta$ is equivalent to the category of right 
$k\Delta/(R)$-modules we deduce $\Mod(\Delta)\cong \Mod(D)$.

Let $\Mscr \in \Xscr_1$ and $M= 
\Ext^{1}_{\P2q}(\Escr,\Mscr)$. 
By functoriality, multiplication by $x,y,z \in A$ induces 
linear maps $M(\lambda_{-1}) : \Ext^{1}(\Oscr(1), \Mscr) \r 
\Ext^{1}(\Oscr, \Mscr)$ and 
$M(\lambda_{-2}) : \Ext^{1}(\Oscr(2), \Mscr) \r \Ext^{1}(\Oscr(1), 
\Mscr)$ ($\lambda = X,Y,Z$). 
Hence $M$ is determined by the following representation of $\Delta$ 
\begin{eqnarray*} 
\begin{array}{ccccc} 
\quad & {\stackrel{M(X_{-2})}{\longrightarrow}} & \quad & 
{\stackrel{M(X_{-1})}{\longrightarrow}} & 
\quad \\ 
H^1(\P2q,\Mscr(-2)) & 
{\stackrel{M(Y_{-2})}{\longrightarrow}} & 
H^1(\P2q,\Mscr(-1)) & 
{\stackrel{M(Y_{-1})}{\longrightarrow}} &  H^1(\P2q,\Mscr) \\ 
\quad & {\stackrel{M(Z_{-2})}{\longrightarrow}} & \quad & 
{\stackrel{M(Z_{-1})}{\longrightarrow}} & 
\quad 
\end{array} 
\end{eqnarray*} 
For further use we note that the Euler form $\chi( S_i,S_j)$ 
is given by the following matrix 
\begin{equation} 
\label{ref-5.4-44} 
\begin{pmatrix} 
1 & -3 & 3\\ 
0 &1& -3\\ 
0 & 0 & 1 
\end{pmatrix} 
\end{equation} 
where $i$ refers to the columns and $j$ refers to the rows. 

Let $\Delta^0$ be the full subquiver of $\Delta$ consisting of the 
vertices $-2,-1$ and let $\Res:\Mod(\Delta)\r \Mod(\Delta^0)$ be the 
obvious restriction functor. $\Res$ has a left adjoint 
which we denote 
by $\Ind$. If $e$ is the sum of the vertices of $\Delta^0$ then 
$\Ind=-\otimes_{k\Delta^0} e k\Delta$. Note that $\Res\circ \Ind =\Id$. 

The following  was already observed by Le Bruyn in the case of the 
homogenized Weyl algebra. 
\begin{lemmas} \label{ref-5.1.1-45} 
If $\Mscr\neq \Oscr$ is a normalized
line bundle on $\P2q$ and 
$M=\Ext^1(\Escr,\Mscr)$ then $M=\Ind\Res M$. 
\end{lemmas} 
\begin{proof} 
This follows from an 
argument by Baer \cite[Corollary 7.2]{Baer}. For the convenience of the 
reader we repeat this argument. 

We  say that two objects $A$, $B$ in an abelian category are 
orthogonal ($A\perp B$) if $\Hom(A,B)=\Ext^1(A,B)=0$. 

We have $\RHom(\Escr,\Mscr)=M[-1]$, $\RHom(\Escr,\Oscr)=S_0$. Thus 
$\Ext^i(\Mscr,\Oscr)=\Ext^i(M[-1],S_0)=\Ext^{i+1}(M,S_0)$. In 
particular $\Hom(M,S_0)=0$ and \linebreak
$\Ext^1(M,S_0)=\Hom(\Mscr,\Oscr)=H^2(\P2q,\Mscr(-3))^\ast=0$ where we 
have used Serre duality and Theorem \ref{ref-3.5-31}.  We conclude by 
lemma \ref{ref-5.1.2-46} below. 
\end{proof} 
\begin{lemmas} \label{ref-5.1.2-46} 
 Let $M\in\mod \Delta$. Then $M=\Ind \Res M$ if and only 
if $M\perp S_0$. 
\end{lemmas} 
\begin{proof} 
First assume $M=\Ind \Res M$. Put $M^0=\Res M$ and take a projective 
resolution 
\[ 
0\r F^0_1\r F^0_0\r M^0\r 0 
\] 
Applying $\Ind$ we get a projective resolution of $M$ of the form 
\[ 
0\r S_0^a\r F_1\r F_0 \r M\r 0 
\] 
for some $a\in\NN$ where $F_i=\Ind F_i^0$. The fact that 
$\Hom(F_1,S_0)=\Hom(F_0,S_0)= 0 $ (by adjointness)  implies 
$\Hom(M,S_0)=0$ and $\Ext^1(M,S_0)=0$. 

To prove the converse let $N=\Ind \Res M$. By adjointness we have a map 
$p:N\r M$ whose kernel 
$K$ and cokernel $C$  are direct sums of $S_0$. We have 
$\Hom(M,S_0)=0$ and hence $\Hom(C,S_0)=0$. Thus $C=0$ and $p$ is surjective. 

Applying $\Hom(-,S_0)$ to the short exact sequence 
\[ 
0\r K \r N\r M \r 0 
\] 
and using $\Hom(N,S_0)=0$ (by adjointness) yields $\Hom(K,S_0)=0$ and 
hence $K=0$. Thus $p$ is an isomorphism and we are done. 
\end{proof}

\begin{lemmas} 
\label{ref-5.1.3-47} 
Let $p = (\alpha, \beta, \gamma)\in E$ and put $(\alpha_{i}, 
\beta_{i}, 
\gamma_{i}) = p^{\sigma^{i}}$. $p$ corresponds to a point module $P$ 
of $A$. Put $\Pscr=\pi P$ and $\bar{P}=\omega\Pscr$. 
\begin{enumerate} 
\item  $H^i(\P2q,\Pscr(n))=0$ for all $n$ and $i>0$. In particular $\Pscr\in 
\Xscr_0$. 
\item $\dim\bar{P}_n=1$ for all $n$ and $\bar{P}_{\ge n}$ is a shifted 
point module for all $n$. In particular $\bar{P}_{\ge 0}=P$. 
\item $H^0(\P2q,\Pscr(n))=\bar{P}_n$. 
\item The representation of $\Delta$ corresponding to $\Pscr$  is 
\begin{eqnarray*} 
\begin{array}{ccccc} 
\quad & {\stackrel{\alpha_{-2}}{\longrightarrow}} & \quad & 
{\stackrel{\alpha_{-1}}{\longrightarrow}} & 
\quad \\ 
k & {\stackrel{\beta_{-2}}{\longrightarrow}} & 
k & 
{\stackrel{\beta_{-1}}{\longrightarrow}} &  k \\ 
\quad & {\stackrel{\gamma_{-2}}{\longrightarrow}} & \quad & 
{\stackrel{\gamma_{-1}}{\longrightarrow}} & 
\quad 
\end{array} 
\end{eqnarray*} 
\item Denote the representation in the previous diagram also by 
$p$. We have $p=\Ind \Res p$. 
\end{enumerate} 
\end{lemmas} 
\begin{proof} 
\begin{enumerate} 
\item Since the $\Pscr(n)$ are all obtained from point modules, it 
suffices to  treat the case $n=0$. We use lemma \ref{ref-2.8.3-19} and the 
discussion before that. We have $\Pscr=i_\ast \Oscr_p$ 
and hence $\Ext^j(\Oscr,\Pscr)=\Ext^j_E(Li^\ast\Oscr,\Oscr_p)= 
\Ext^j_E(\Oscr_E,\Oscr_p)=0$ for $j>0$. 
\item This is easy to check. 
\item 
Use $\omega = \underline{\Hom}_{\Tails}(\Oscr,-)$.
\item This follows from the previous step. 
\item According to lemma \ref{ref-5.1.2-46} we need $\Ext^i(p,S_0)=0$ for 
$i\le 1$. This follows from the fact that we have 
$\Ext^i(p,S_0)=\Ext^i(\Pscr,\Oscr)=0$ for $i\le 1$ by lemma 
\ref{ref-3.4-29}.\qed 
\end{enumerate} 
\def\qed{}\end{proof} 

To simplify the discussion below we define $\Rscr_n$ (for $n\ge 1$) as 
the category in 
which the objects are the 
normalized line bundles on 
$\P2q$ with invariant $n$ 
and the morphisms are the isomorphisms in $\coh(\P2q)$. Thus 
$\Rscr_n$ is  a groupoid. Note that we do not know yet if $\Rscr_n\neq 
\emptyset$. This question will be addressed below. 


\subsection{$\Rscr_n$ is non-empty} 
\label{ref-5.2-48} 
To simplify things we assume that $\sigma$ has infinite 
order. The following result is necessary for the dimension 
computations in \S\ref{ref-5.5-53}. 
\begin{lemmas} \label{ref-5.2.1-49} 
The set $\Rscr_n$ is not empty. 
\end{lemmas} 
\begin{proof} Let $\Sscr$ be a line object on $\P2q$. Writing 
$\Sscr$ as the cokernel of a map $\Oscr(-1)\r \Oscr$ we find by 
Theorem \ref{ref-2.6.2-13} that if $n\ge -1$ then $H^0(\P2q,\Sscr(n))$ 
has dimension $n+1$. 

By \cite{ATV2} there exist at most three line objects 
$\Sscr'$ such that $\Sscr'(-1)$ is a subobject of $\Sscr$ and 
furthermore these three line objects contain in turn any other object 
contained in $\Sscr$. 

Hence if $n\ge 0$ then we may pick 
an epimorphism $f:\Oscr\r \Sscr(n)$ (a generic $f$ will do). 
Put $\Iscr=(\ker f)(1)$. Using Proposition \ref{ref-3.1-24} we 
find $[\Iscr(-1)]=[\Oscr]-([\Sscr]+n[\Pscr])$ and hence 
$[\Iscr]=[\Oscr]-n[\Pscr]$.  It is easy to see that $\Iscr$ is 
reflexive.  
Thus $\Iscr\in \Rscr_n$. 
\end{proof} 
Below we will show that $\Rscr_n$ is parametrized by an algebraic 
variety of dimension $2n$. The amount of freedom in the construction 
exhibited in the proof of lemma \ref{ref-5.2.1-49} is less than or equal 
to $2(\text{choice of $\Sscr$})+n(\text{choice of $f$})$ parameters, hence for 
$n>2$ this construction can not possibly yield all elements of 
$\Rscr_n$. In 
\S\ref{ref-5.6-66} we will 
 exhibit a related construction which works for all $n$. 
\subsection{First description of $\Rscr_n$} 
Let $\Cscr_n$ be the image of $\Rscr_n$ under the equivalence 
$\Xscr_1\cong\Yscr_1$. 

\begin{theorems} \label{ref-5.3.1-50} Let $n\ge 1$. 
There is  an equivalence of categories 
\begin{eqnarray} \label{ref-5.5-51} 
\Rscr_n {{{\Ext_{\P2q}^{1}(\Escr, -)} \atop 
{\lr}} \atop {{\longleftarrow} \atop {\Tor_{1}^{{{D}}}(-, \Escr)}}} 
\Cscr_n 
\end{eqnarray} 
where 
\begin{multline*} 
\Cscr_n = \{ M \in \mod(\Delta) \mid  \underline{\dim} M=(n,n,n-1) 
\text{ and } \\ 
\Hom_{\Delta}(M,p) = 0, \Hom_{\Delta}(p,M) = 0 \mbox{ for all } p \in E\}. 
\end{multline*} 
\end{theorems} 
\begin{proof} 
First, let $\Mscr$ be an object of $\Rscr_n$. By Proposition 
\ref{ref-4.3-39} we have $i^{*}\Mscr \cong \Nscr$ for a 
line bundle $\Nscr$ of degree zero on $E$. Hence (for all $p \in E$) we 
have $\RHom_{E}(Li^{*}\Mscr, \Oscr_{p}) = k$. Since 
\[ 
\RHom_{E}(Li^{*}\Mscr, \Oscr_{p}) \cong \RHom_{\P2q}(\Mscr, 
i_{*}\Oscr_{p}) 
\] 
we obtain $\RHom_{\P2q}(\Mscr, \Pscr) = k$ 
where $\Pscr = i_{*}\Oscr_{p}$ is the corresponding point object
on $\P2q$. Writing $M = \Ext^{1}_{\P2q}(\Escr, \Mscr)$ this means 
that $\RHom_{{{D}}}(M[-1], p) = k$, proving that $\Hom_{{{D}}}(M,p) = 0$ and 
$\Ext^2_{{D}}(M,p)=0$. 
By BKS-duality (Theorem \ref{ref-2.9.1-23}) we obtain $\Hom_{{D}}(p',M)=0$ for 
some other point $p'$ determined by $p$ (and determining $p$). Hence 
$M\in\Cscr_n$. 

Conversely, let $M$ be an object of $\Cscr_n$.  Thus (using Serre 
duality on $\P2q$ again) $\Hom_{{D}}(M,p)=\Ext^2_{{D}}(M,p)=0$ for all 
$p\in E$. Now $\gldim 
{{D}}=2$ so we may 
compute $\dim \Ext^1_{{D}}(M,p)$ using the Euler form 
\eqref{ref-5.4-44} on $\mod({{D}})$. We 
obtain \linebreak
$\Ext^1_{{D}}(M,p)=k$.  
In other words $\RHom_{{{D}}}(M[-1],p) = k$. 

Put $\Mscr = M[-1] \LotimesD \Escr$.  By the category 
equivalence between $D^b(\coh(\P2q))$ and $D^b(\mod({{D}}))$ we obtain 
$\RHom_{\P2q}(\Mscr, \Pscr) = k$, giving (by adjointness) \linebreak
$\RHom_{E}(Li^{*}\Mscr, \Oscr_{p}) = k$. Since $E$ is a smooth 
elliptic curve it is easy to see that this implies that $Li^\ast 
\Mscr$ is a line bundle. Hence by Propositions \ref{ref-4.3-39} and 
\ref{ref-4.4-40} the same is true for $\Mscr$. 
\end{proof} 
\subsection{Application} 
Using the material in the previous sections it is now easy to parametrize 
the line bundles on $\P2q$ with invariant one. 
\begin{theorems}  The representations in $\Cscr_1$ are the 
representations 
\begin{eqnarray} 
\label{ref-5.6-52} 
\begin{array}{ccccc} 
\quad & {\stackrel{\alpha}{\longrightarrow}} & \quad & 
{\stackrel{0}{\longrightarrow}} & 
\quad \\ 
k & {\stackrel{\beta}{\longrightarrow}} & 
k & 
{\stackrel{0}{\longrightarrow}} &  0 \\ 
\quad & {\stackrel{\gamma}{\longrightarrow}} & \quad & 
{\stackrel{0}{\longrightarrow}} & 
\quad 
\end{array} 
\end{eqnarray} 
for some $(\alpha,\beta,\gamma)\in \PP^2-E$ 
\end{theorems} 
\begin{proof} 
First let  $F\in 
\Cscr_1$. $F$ is given by a representation as in \eqref{ref-5.6-52}. 
Then the condition 
$\Hom_\Delta(p,F)=0$ for $p\in E$ implies 
$(\alpha,\beta,\gamma)\not \in E$. 

Conversely let $F$ be as in \eqref{ref-5.6-52} with 
$(\alpha,\beta,\gamma)\not\in E$. Then we immediately have 
$\Hom_\Delta(p,F)=\Hom_\Delta(F,p)=0$ for $p\in E$. 
\end{proof} 

\subsection{Second description of $\Rscr_n$} 
\label{ref-5.5-53} 
Although the category $\Cscr_n$ has a fairly elementary description, it 
is not so easy to handle.  In particular the analogy with the Weyl 
algebra case is not obvious. 
We will now give another description of 
$\Rscr$ which is more similar to the one used for the Weyl algebra. In 
particular it will follow that the isomorphism classes of objects in 
$\Rscr_{n}$ are parametrized by smooth affine varieties of dimension 
$2n$. 

In general, let $Q$ be a quiver without oriented cycles and write 
$Q_{0}$, $Q_1$ for respectively the set of 
vertices and edges of $Q$. Let ``$\cdot$'' be the standard scalar product on 
$\ZZ^{Q_0}$: $(\alpha_v)_v\cdot (\beta_v)_v=\sum_v \alpha_v\beta_v$. 

Let $\theta \in \ZZ^{Q_{0}}$. 
A representation $F$ of $Q$ is called 
$\theta$-{\it semistable} (resp.\ {\it stable}) if $\theta\cdot 
\underline{\dim} F=0$ and $ \theta\cdot 
\underline{\dim} N\geq 0$ (resp. $> 0$) for every 
non-trivial subrepresentation $N$ of $F$. 

We have $K_0(\mod 
Q)=\ZZ^{Q_0}$, canonically. 
It is a fundamental fact \cite{schofield} 
that $F$ is semistable for some $\theta$ if and only there exists 
$G\in\mod(Q)$ such that $F\perp G$. The relation 
between $\theta$ and $\underline{\dim} G$ is such that the forms 
$-\cdot \theta$ and $\chi(-,\underline{\dim} G)$ are proportional. 

Fix a dimension vector $\alpha\in \ZZ^{Q_0}$ and let $\Rep(Q,\alpha)$ 
be the corresponding representation space, 
i.e. $\Rep(Q,\alpha)=\prod_{i\in Q_1} M_{\alpha_{h(i)}\times 
\alpha_{t(i)}}(k)$ where the maps $h,t:Q_1\r Q_0$ associate to an arrow 
its  begin and 
end vertex.  The isomorphism class of representations of 
dimension vector $\alpha$ are in one-one correspondence with the 
orbits of the group $\Gl(\alpha)=\prod_{v\in Q_0} 
\Gl_{\alpha_v}(k)$ acting on $\Rep(Q,\alpha)$ by conjugation. 

Associated to $G\in \mod(Q)$ there is a semi-invariant function 
$\phi_G$ on $\Rep(Q,\alpha)$ such that the set 
\begin{equation} 
\label{ref-5.7-54} 
{}^\perp G=\{F\in\Rep(Q,\alpha)\mid F\perp G\} 
\end{equation} 
coincides with $\{\phi_G\neq 0\}$. In particular \eqref{ref-5.7-54} is 
affine. 

\begin{lemmas}  \label{ref-5.5.1-55} 
There exists $V\in \mod(\Delta^0)$ with 
$\underline{\dim} V=(6,3)$ 
such that 
\begin{enumerate} 
\item 
for all $M\in \Cscr_n$  we have 
$M^0\perp V$ 
where $M^0=\Res M$ and 
\item if $p\in E$ then $\RHom_{\Delta^0}(\Res p,V)\neq 0$. 
\end{enumerate} 
\end{lemmas} 
\begin{proof} 
\begin{enumerate} 
\item 
Pick a degree zero line bundle $\Uscr$ on $E$ which is not of the form 
$\Oscr((o)-(3n\xi)) $ for $n\in \NN$ (where $o$, $\xi$ are as in 
Proposition \ref{ref-4.3-39}). 

Let $\Mscr\in \Rscr_n$. Then we have by adjointness 
$\RHom(\Mscr,i_\ast \Uscr)= \linebreak 
\RHom_E(Li^\ast \Mscr,\Uscr)$. By Proposition 
\ref{ref-4.3-39} we have $Li^\ast \Mscr= \Oscr((o)-(3n\xi)) $. We 
conclude by Serre duality for $E$ that $\RHom(\Mscr,i_\ast\Uscr)=0$. 
Now put $M=\Ext^1(\Escr,\Mscr)$ and $U'=\RHom(\Escr,i_\ast\Uscr)$.  
We obtain \linebreak
$\RHom_{{D}}(M[-1],U')=0$. 

What is $U'$? By adjointness we have 
$\RHom(\Escr,i_\ast\Uscr)=\RHom_E(Li^\ast\Escr,\Uscr)$. An easy 
verification shows that 
$Li^\ast\Escr=\sigma_\ast^2(\Lscr)\otimes\sigma_\ast(\Lscr)\oplus 
\sigma_\ast\Lscr\oplus \Oscr_E$. Thus by Riemann-Roch and Serre 
duality $U'=U[-1]$ where $\underline{\dim} U=(6,3,0)$.  Put $V=\Res 
U$. Thus $\underline{\dim}V=(6,3)$. 

Replacing $M$ with a projective resolution it is easy to see that \linebreak
$\RHom_\Delta(M,U)=\RHom_{\Delta^0}(M^0,V)$.  It follows that 
$\Hom_{\Delta^0}(M^0,V)=0$ and $\Ext^1_{\Delta^0}(M^0,V)=0$. 
\item Put $Q=\Res 
p$ for $p\in E$. Then $\RHom_{\Delta^0}(Q,V)=\RHom_{\Delta}(p, 
U)=\RHom_{\Delta}(p[-1], U')= \RHom_{\P2q}(i_\ast\Oscr_p[-1], i_\ast 
\Uscr)=\RHom_E(Li^\ast i_\ast \Oscr_p[-1],\Uscr)$. Now $Li^\ast 
i_\ast \Oscr_p[-1]$ is a nonzero complex whose homology has finite 
length. It is easy to deduce from this $\RHom_E(Li^\ast i_\ast 
\Oscr_p[-1],\Uscr)\neq 0$. Hence we are done \qed 
\end{enumerate} 
\def\qed{}\end{proof} 
We obtain the following consequence. 
\begin{lemmas} \label{ref-5.5.2-56} If $M\in \Cscr_n$ and $M^0=\Res M$ then 
$M^0$ is 
$\theta$-semistable for $\theta=(-1,1)$. 
\end{lemmas} 
\begin{proof} This is a straightforward verification. 
\end{proof} 
\begin{lemmas} \label{ref-5.5.3-57} Assume that $\sigma$ has infinite 
order. Let $N$ be a 
representation of $\Delta^0$ of dimension vector $(n,n)$, $n\ge 1$. If 
$\Hom_{\Delta^0}(N,\Res 
p)=\Hom_{\Delta^0}(\Res p,N)=0$ for all $p\in E$ then $\dim (\Ind 
N)_0\le n-1$. 
\end{lemmas} 
\begin{proof} Assume the lemma is false. Thus $\dim (\Ind N)_0\ge 
n$. Then we may construct a surjective map $\Ind N\r W$ where 
$\underline{\dim} W=(n,n,n)$. We will consider $W\LotimesD \Escr$ 
and $Li^\ast(W\LotimesD \Escr)$. Note that since $E$ is smooth, 
$Li^\ast(W\LotimesD 
\Escr)$ is the sum of its homology. 

We have for $p\in E$: 
\begin{equation} 
\label{ref-5.8-58} 
\begin{split} 
\Ext^j_E(Li^\ast(W\LotimesD \Escr), \Oscr_p) 
&=\Ext^j_{\P2q}(W\LotimesD \Escr,i_\ast \Oscr_p)\\ 
&=\Ext^j_\Delta(W,p) 
\end{split} 
\end{equation} 
Now a simple computation shows that $\chi(W,p)=0$.  Furthermore we 
have \linebreak
$\Hom_\Delta(W,p)\subset \Hom_\Delta(\Ind 
N,p)=\Hom_{\Delta^0}(N,\Res p)=0$. Finally by Serre duality on $\P2q$ (see 
Theorem \ref{ref-2.9.1-23}) we have 
$\Ext^2_\Delta(W,p)=\Hom_{\Delta}(p',W)^\ast=0$. We conclude that also 
$\Ext^1_\Delta(W,p)=0$.  It follows from \eqref{ref-5.8-58} that 
$Li^\ast(W\LotimesD \Escr)=0$.  Hence by lemmas \ref{ref-2.8.4-20} and 
\ref{ref-2.8.2-18} we deduce $W\LotimesD \Escr=0$ and hence $W=0$ which is 
a contradiction. 
\end{proof} 
We can now prove our main result. 
\begin{theorems} 
\label{ref-5.5.4-59} 
Assume that $\sigma$ has infinite order. Let $V\in 
\mod(\Delta^0)$ be as in lemma \ref{ref-5.5.1-55}. 
\begin{enumerate} 
\item 
The functors $\Res$ and $\Ind$ define inverse equivalences between 
$\Cscr_n$ and the 
following category 
\[ 
\Dscr_n=\{F\in \mod (\Delta^0)\mid \underline{\dim} F=(n,n), 
F\perp V, \dim (\Ind F)_0\ge n-1 \} 
\] 
\item 
The representations in $\Dscr_n$ are $\theta$-stable for 
$\theta=(-1,1)$. 
\end{enumerate} 
\end{theorems} 
\begin{proof} Below we use often implicitly the already proved 
equivalence $\Cscr_n\cong \Rscr_n$ \linebreak
(Theorem \ref{ref-5.3.1-50}). 
\begin{step} \label{ref-1-60}  $\Res(\Cscr_n)\subset \Dscr_n$. This follows 
from lemmas 
\ref{ref-5.1.1-45} and \ref{ref-5.5.1-55}. 
\end{step} 
\begin{step} \label{ref-2-61} Let $F\in \Dscr_n$. If 
$F'\subset F$ is such that $\underline{\dim} F'=(m,m)$ then 
$\Hom_{\Delta^0}(F',V)=\Ext^1_{\Delta^0}(F',V) 
=\Hom_{\Delta^0}(F/F',V)=\Ext^1_{\Delta^0}(F/F',V)=0$. This follows 
easily by using the Euler form. 
\end{step} 
\begin{step} \label{ref-3-62} Let $F\in \Dscr_n$ and $p\in E$. Then 
$\Hom_{\Delta^0}(F,\Res 
p)=\Hom_{\Delta^0}(\Res p, F)=0$.  Both claims are similar so we only 
consider the first one. 
We have $F\perp V$ hence
$F$ is $\theta$-semistable 
and $\Res p$ is obviously stable. So if $\Hom_{\Delta^0}(F,\Res 
p)\neq 0$ then there is an epimorphism $F\r \Res p$. By Step 
\ref{ref-2-61}  we obtain $\RHom_{\Delta^0}(\Res p,V)=0$. But 
this contradicts the choice of $V$, finishing the argument. 
\end{step} 
\begin{step} \label{ref-4-63} $\Ind(\Dscr_n)\subset \Cscr_n$. Let $F\in 
\Dscr_n$. 
By Step \ref{ref-3-62} and lemma \ref{ref-5.5.3-57} we obtain $\dim (\Ind 
F)_0=n-1$. \\
It remains to show that for $p\in E$ we have $\Hom_\Delta(\Ind F, 
p)=\Hom_\Delta(p,\Ind F)=0$. By lemma \ref{ref-5.1.3-47} we have 
$p=\Ind \Res p$. Thus  $\Hom_\Delta(\Ind F, 
p)=\Hom_{\Delta^0} (F, \Res p)=0$ and similarly $\Hom_\Delta(p,\Ind 
F)=\Hom_{\Delta^0}(\Res p, \Res \Ind F)=\Hom_{\Delta^0}(\Res p, 
F)=0$ where we have used Step \ref{ref-3-62} again. 
\end{step} 
\begin{step}  $\Ind$ and $\Res$ are inverses to each other.  To prove 
this we only need to show $\Ind \circ \Res(F)=F$ for $F\in 
\Cscr_n$. This follows from lemma \ref{ref-5.1.1-45}. 
\end{step} 
\begin{step} Let $F\in \Dscr_n$. Then $F$ is $\theta$-stable. Put a 
filtration $0=F_0\subsetneq 
F_1\subsetneq\cdots \subsetneq F_{m-1}\subsetneq F_m=F$ on $F$ such 
that $F_{i}/F_{i-1}$ is $\theta$-stable. With the same proof as Step 
\ref{ref-3-62} 
it follows that $\Hom_{\Delta^0}(F_{i}/F_{i-1},\Res 
p)=\Hom_{\Delta^0}(\Res p, F_{i}/F_{i-1})=0$ for $p\in E$. Assume 
$\underline{\dim}F_{i}/F_{i-1}=(d_i,d_i)$. Then by lemma \ref{ref-5.5.3-57} 
we have $\dim (\Ind(F_{i}/F_{i-1}))_{0}\le d_i-1$.  From the right 
exactness of $\Ind$ we deduce $\dim (\Ind F)_{0}\le n-m$. Hence $m=1$ and 
thus $F$ is stable. \qed 
\end{step} 
\def\qed{}\end{proof} 
Below we will define some varieties. We take the classical 
viewpoint. So they are always reduced. 

Let $V$ as in lemma \ref{ref-5.5.1-55}. Let $\alpha=(n,n)$ and put 
\begin{equation} 
\label{ref-5.9-64} 
\begin{split}
\tilde{D}_n & =\{F\in \Rep(\Delta^0,\alpha)\mid F\in 
\Dscr_n\} \\
& =\{F\in\Rep(\Delta^0,\alpha)\mid \phi_V(F)\neq 0, \dim(\Ind 
F)_0\ge n-1\}. 
\end{split}
\end{equation} 
It is clear that $\tilde{D}_n$ is a closed subset of $\{\phi_V\neq 
0\}$ so in particular $\tilde{D}_n$ is affine. Put 
$D_n=\tilde{D}_n\quot \Gl(\alpha)$. 
\begin{theorems} 
\label{ref-5.5.5-65} The affine variety  $D_n$ is smooth of dimension $2n$. 
The isomorphism classes in $\Dscr_n$ (and hence in 
$\Cscr_n$ and $\Rscr_n$) are in natural bijection with the points in 
$D_n$. 
\end{theorems} 
\begin{proof}  Since all representions in 
$\tilde{D}_n$ are stable by Theorem \ref{ref-5.5.4-59}, all 
$\Gl(\alpha)$-orbits on $\tilde{D}_n$ are 
closed and so $D_n$ is really the orbit space for the $\Gl(\alpha)$ 
action on $\tilde{D}_n$. This proves that the isomorphism classes in 
$\Dscr_n$ are in natural bijection with the points in $D_n$. 

To prove that $D_n$ is smooth it suffices to prove that $\tilde{D}_n$ 
is smooth 
(this follows for example using the Luna slice theorem \cite{Luna}). 

We first estimate the dimension 
of $\tilde{D}_n$.  We write the equations of $A$ in the usual form 
$M(xyz)^t$. Given $n\times n$-matrices $X$, $Y$, $Z$ 
let $M(X,Y,Z)$ be obtained from $M$ by replacing $(x,y,z)$ by $X,Y,Z$ 
(thus $M(X,Y,Z)$ is a $3n\times 3n$-matrix). Then $\tilde{D}_n$ has 
the following alternative description: 
\[ 
\tilde{D}_n=\{ (X,Y,Z)\in M_n(k)^3\mid \phi_V(X,Y,Z)\neq 0\text{ and }\rk 
M(X,Y,Z)\le 3n-(n-1)\}. 
\] 
By \S\ref{ref-5.2-48} $\tilde{D}_n$ is non-empty.  The triples satifying 
$\phi_V(X,Y,Z)\neq 0$ are a dense open subset of $M_n(k)^3$ and 
hence they represent a variety of dimension $3n^2$. Imposing that 
$M(X,Y,Z)$ should have corank $\ge n-1$ represents $(n-1)^2$ 
independent conditions. So the irreducible components of $ \tilde{D}_n 
$ have dimension $\ge 
3n^2-(n-1)^2$. 

Define $\tilde{C}_n$ by 
\[ 
\{G\in \Rep(\Delta,\tilde{\alpha})\mid  G\cong \Ind \Res G, \Res G\in 
\tilde{D}_n\} 
\] 
where $\tilde{\alpha}=(n,n,n-1)$ (as usual we assume the points of 
$\Rep(\Delta,\tilde{\alpha})$ to satisfy the relation imposed on $\Delta$). 

To extend $F\in \tilde{D}_n$ to a point in $\tilde{C}_n$ we need to 
choose a basis in $(\Ind F)_0$. Thus $\tilde{C}_n$ is a principal 
$\Gl_{n-1}(k)$ 
fiber bundle over $\tilde{D}_n$. In particular $\tilde{C}_n$ is smooth 
if and only $\tilde{D}_n$ is smooth and the irreducible components of 
$\tilde{C}_n$ have dimension $\ge 
3n^2-(n-1)^2+(n-1)^2=3n^2$. Note that by the description of $\Cscr_n$ 
in Theorem \ref{ref-5.3.1-50} it follows that $\tilde{C}_n$ is an open subset of 
$\Rep(\Delta,\tilde{\alpha})$. 

Let $x\in \tilde{C}_n$. The stabilizer of $x$ consists of scalars thus 
if we put ${{G}}=\Gl(\tilde{\alpha})/k^\ast$ then we have inclusions 
$\Lie(G)\subset T_x(\tilde{C}_n)= 
T_x(\Rep(\Delta,\tilde{\alpha}))$. Voigt in \cite[Ch. 2, \S3.4]{voigt} has 
shown that there 
is a natural inclusion 
$T_x(\Rep(\Delta,\tilde{\alpha}))/\Lie(G)\hookrightarrow 
\Ext^1_\Delta(x,x)$ (Voigt actually obtains an isomorphism since he is 
not assuming his representation spaces to be reduced). Now $x$ corresponds 
to some line bundle $\Hscr$ 
on  $\P2q$ and we have $\Ext^1_\Delta(x,x)=\Ext^1(\Hscr,\Hscr)$. An 
easy computation shows $\chi(\Hscr,\Hscr)=\chi(x,x)=1-2n$. We have 
$\Hom(\Hscr,\Hscr)=k$ and by Serre duality 
$\Ext^2(\Hscr,\Hscr)=\Hom(\Hscr,\Hscr(-3))=0$. Thus $\dim 
\Ext^1(\Hscr,\Hscr)=2n$. 

Hence we obtain $3n^2\le \dim T_x(\tilde{C}_n)\le 2n+\dim 
G=2n+2n^2+(n-1)^2-1=3n^2$. Thus $\dim T_x(\tilde{C}_n)=3n^2$ is 
constant and hence $\tilde{C}_n$ is smooth. We also obtain $\dim 
\tilde{D}_n=3n^2-(n-1)^2$. 

The dimension of $D_n$ is equal to $\dim \tilde{D}_n-\dim 
\Gl(\alpha)+1=3n^2-(n-1)^2 -2n^2+1=2n$. This finishes the proof. 
\end{proof} 

\subsection{Explicit construction of the elements in $\Rscr_n$} 
\label{ref-5.6-66} 
For simplicity we assume throughout that $\sigma$ has infinite order. 

In the discussion below we have to compute the cohomology of a line 
object. 
\begin{lemmas} Let $\Sscr=\pi(A/uA)$ be a line object on $\P2q$. Let $m\le 
-1$. Then $H^1(\P2q,\Sscr(m))\cong (A/Au)_{-m-2}^\ast$ and 
$H^i(\P2q,\Sscr(m))=0$ for $i\neq 1$. Furthermore if 
$\eta\in A_1$ then the induced linear map 
$H^1(\P2q,\Sscr(m))\xrightarrow{\cdot\eta} 
H^1(\P2q,\Sscr(m+1))$  corresponds to $(\eta\cdot)^\ast$ on $(A/Au)^\ast$. 
\end{lemmas} 
\begin{proof} 
That $H^0(\P2q,\Sscr(m))=0$ follows by writing $\Sscr$ as the cokernel 
of a map $\Oscr(-1)\r \Oscr$ and invoking Theorem 
\ref{ref-2.6.2-13}. That $H^2(\P2q,\Sscr(m))=0$ follows by 
Serre duality (Theorem \ref{ref-2.9.1-23}). 

Using Theorem \ref{ref-2.6.2-13} we find $H^1(\P2q,\Sscr(m))=\ker ( 
\Ext^2(\Oscr(-m),\Oscr(-1))\xrightarrow{(-,u\cdot)} 
\Ext^2(\Oscr(-m),\Oscr))$. Using Serre duality (Theorem \ref{ref-2.9.1-23}) 
this translates into \linebreak
$H^1(\P2q,\Sscr(m))=\coker(\Hom(\Oscr(-1),\Oscr(-m-3))^\ast 
\xrightarrow{(u\cdot,-)^\ast} 
\Hom(\Oscr,\Oscr(-m-3))^\ast)$. Dualizing yields that indeed 
$H^1(\P2q,\Sscr(m))\cong  (A/Au)_{-m-2}^\ast$. That $\eta$ acts in the indicated way 
follows by inspecting the appropriate commutative diagram. 
\end{proof} 
\begin{corollarys} Let $\Sscr=\pi(A/uA)$ be a line object on 
$\P2q$. Then $\Sscr(-1)$ corresponds to $S[-1]$ where $S$ is the 
representation of $\Delta$ given by 
\begin{eqnarray} 
\label{ref-5.10-67} 
\begin{array}{ccccc} 
\quad & {\stackrel{(x\cdot)^\ast}{\longrightarrow}} & \quad & 
{\stackrel{}{\longrightarrow}} & 
\quad \\ 
(A/Au)_{1}^\ast & 
{\stackrel{(y\cdot)^\ast}{\longrightarrow}} & 
k & 
{\stackrel{}{\longrightarrow}} &  0 \\ 
\quad & {\stackrel{(z\cdot)^\ast}{\longrightarrow}} & \quad & 
{\stackrel{}{\longrightarrow}} & 
\quad 
\end{array} 
\end{eqnarray} 
\end{corollarys} 
Since line objects on $\P2q$ are of the form $\pi(A/uA)$ they 
are naturally parametrized by points in $\PP(A_1)$. 
\begin{propositions} \label{ref-5.6.3-68} Let $\Iscr$ be a normalized line 
bundle on $\P2q$ with 
invariant $n>0$. Then the set of line objects $\Sscr$ such that 
$\Hom(\Iscr,\Sscr(-1))\neq 0$ is a curve of degree $n$ in 
$\PP(A_1)$. In particular this set is non-empty. 
\end{propositions} 
\begin{proof} Let $\Sscr=\pi(A/uA)$ with $u=\alpha x+\beta y+\gamma 
z$. 
Put $I=\Ext^1(\Escr,\Iscr)$, \linebreak
$S=\Ext^1(\Escr,\Sscr(-1))$. Then 
$\Hom(\Iscr,\Sscr(-1))=\Hom_\Delta(I,S)=\Hom_{\Delta^0}(\Ind 
I^0,S)=\Hom_{\Delta^0}(I^0,S^0)$ where 
$I^0$, $S^0$ are the restrictions of $I$ and $S$ to $\Delta^0$. 

Assume that $I^0$ is given by matrices $X,Y,Z\in M_n(k)$. Then an easy 
verification shows that $\Hom_{\Delta^0}(I^0,S^0)\neq 0$ if and only if 
$\det(\alpha X+\beta Y+\gamma Z)=0$. This is a homogeneous equation in 
$(\alpha,\beta,\gamma)$ and we have to show that it is not identically 
zero, i.e. we have to show that there is at least one $\Sscr$ such 
that $\Hom(\Iscr,\Sscr(-1))=0$. This follows from lemma \ref{ref-5.6.4-69} 
below. 
\end{proof} 

\begin{lemmas} 
\label{ref-5.6.4-69} 
Let $\Iscr$ be a normalized line bundle on $\P2q$ with invariant 
$n$ and let $\Pscr$ be a point object on $\P2q$. Then, modulo 
zero dimensional objects, there exist at 
most $n$ different line objects $\Sscr$ such that 
$\Hom(\Iscr,\Sscr(-1))\neq 0$ and such that $\Hom(\Sscr,\Pscr)\neq 
0$. 
\end{lemmas} 
\begin{proof}  We use induction on $n$. Writing $\Sscr$ as the 
cokernel of a map $\Oscr(-1)\r \Oscr$ we deduce by Theorem 
\ref{ref-2.6.2-13} that 
$\Hom(\Oscr,\Sscr(-1))=0$. So the case 
$n=0$ is clear by Corollary \ref{ref-3.6-34}. 
Assume $n>0$. Let  $(\Sscr_i)_{i=1,\ldots,m}$ 
be the different line objects (modulo zero dimensional objects) 
satisfying 
$\Hom(\Iscr,\Sscr_i(-1))\neq 0$ and $\Hom(\Sscr_i,\Pscr)\neq 
0$. If $m=0$ then we are done. So assume $m>0$. Let $\Sscr'_i(-1)$ be 
the kernel of a non-trivial map $\Sscr_i\r 
\Pscr$. It is proved in \cite{ATV2} that there is some different 
point object $\Pscr'$ such that for all $i$: 
$\Hom(\Sscr_i,\Pscr')\neq 0$. 
Let $\Iscr'(-1)$ be the kernel of a 
non-trivial map $\Iscr\r \Sscr_1(-1)$.  The subobjects of line 
objects are shifted line objects 
and hence 
the image of $\Iscr$ in $\Sscr_1(-1)$ is a shifted line object.  We 
find by \eqref{ref-3.1-24}: $[\Iscr']=[\Oscr]-(n-b)[\Pscr]$ with 
$b\ge 1$. From this 
we deduce that the invariant of $\Iscr'$ is $\le n-1$. 

Since $\Pscr=i_\ast \Oscr_p$ 
for some point $p\in E$ it follows by adjointness and by Proposition 
\ref{ref-4.3-39} that 
$\dim \Hom(\Iscr,\Pscr(-1))=1$. Hence the composition $\Iscr'(-1)\r\Iscr\r 
\Sscr_i(-1)$ maps $\Iscr'(-1)$ to $\Sscr'_i(-2)$. We claim that for 
$i>1$ this map must be nonzero. If not then there is a non-trivial 
map $\Iscr/\Iscr'(-1)\r  \Sscr_i(-1)$ and since $\Iscr/\Iscr'(-1)$  
is also subobject 
of $\Sscr_1(-1)$ it follows that $\Sscr_1$ and $\Sscr_i$ have a common 
subobject. But this is impossible since $\Sscr_1$ and $\Sscr_i$ are 
different modulo zero dimensional objects. 

Hence $\Hom(\Iscr',\Sscr'_i(-1))\neq 0$ for $i=2,\ldots,m$. Since the 
$\Sscr'_i$ are still different modulo zero dimensional objects, we 
obtain $m-1\le n-1$ and hence $m\le n$. 
\end{proof} 
The following lemma shows how to reduce the invariant of a 
line bundle. 
\begin{lemmas} 
\label{ref-5.6.5-70} Let $\Iscr$ be a normalized line bundle on $\P2q$ 
with invariant $n >0$ . Then there exists a line object $\Sscr$ on $\P2q$ such 
that $\Ext^1(\Sscr(1),\Iscr(-1))\neq 0$. If $\Jscr=\pi J$ is the middle term 
of a corresponding non-trivial extension and $\Jscr^{\ast\ast}=\pi 
J^{\ast\ast}$ 
then $\Jscr^{\ast\ast}$ is a normalized line bundle with invariant 
$\le n-1$. Furthermore $\Jscr^{\ast\ast}/\Iscr(-1)$ is a shifted 
line object. 
\end{lemmas} 
\begin{proof} Using Serre duality we have $\Ext^1(\Sscr(1),\Iscr(-1))= 
\Ext^1(\Iscr(-1),\Sscr(-2))^\ast=\Ext^1(\Iscr,\Sscr(-1))^\ast$. Also using 
Serre duality we deduce  $\Ext^2(\Iscr,\Sscr(-1))=0$. Then a 
simple computation using the Euler form shows that
$\dim\Hom(\Iscr,\Sscr(-1))= 
\dim\Ext^1(\Iscr,\Sscr(-1))$. Hence it follows from Proposition 
\ref{ref-5.6.3-68} that there exist $\Sscr$ such that 
$\Ext^1(\Sscr(1),\Iscr(-1))\neq 0$. 

Now let $\Jscr=\pi J$ be the middle term of a non-trivial extension of 
$\Iscr(-1)$ by $\Sscr(1)$. Then we have 
$[\Jscr]=[\Oscr]-[\Sscr]-n[\Pscr]+[\Sscr]+[\Pscr]=[\Oscr]-(n-1)[\Pscr]$. 

We claim that $\Jscr$ is torsion-free. Assume this is not the case and 
let $\Fscr\subset \Jscr$ be a maximal subobject of $\Jscr$ of 
dimension $\le 1$. So $\Fscr\neq 0$. Since $\Iscr$ is torsion-free we have 
$\Fscr\cap \Iscr(-1)=0$. So we may consider $\Fscr$ as a subobject of 
$\Sscr(1)$. Hence we obtain an extention 
\begin{equation} 
\label{ref-5.11-71} 
0\r \Iscr(-1)\r \Jscr/\Fscr\r \Sscr(1)/\Fscr\r 0 
\end{equation} 
According to lemma \ref{ref-3.4-29} this extension is split. But this 
means that $\Sscr(1)/\Fscr$ is a  subobject of $\Jscr/\Fscr$ of 
dimension $\le 1$, 
contradicting the maximality of $\Fscr$. 

It follows from \cite{ATV2} that $\gkdim 
J^{\ast\ast}/J\le 1$. Thus 
$\Jscr^{\ast\ast}/\Jscr=b[\Pscr]$ for some $b\ge 0$ by Proposition 
\ref{ref-3.2-27}. Hence $[\Jscr^{\ast\ast}]=[\Oscr]-(n-1-b)[\Pscr]$. 

Let $\Sscr'=\Jscr^{\ast\ast}/\Iscr(-1)$. Then by lemma \ref{ref-3.4-29} 
$\Sscr'$ is pure and furthermore we  have $e(\Sscr')=1$. It then 
follows easily using the methods of \cite{ATV2} that $\Sscr'$ is a 
shifted line object. 
This 
finishes the proof. 
\end{proof} 
We can now prove another main result. 
\begin{theorems} 
\label{ref-5.6.6-72} 
Let $\Iscr$ be a normalized line bundle on $\P2q$. Then there exists 
an $m\in\NN$ together with a monomorphism $\Iscr(-m)\hookrightarrow 
\Oscr$ such that there exists a filtration of line bundles $\Oscr=\Mscr_0\supset 
\Mscr_1\supset \cdots 
\supset \Mscr_u=\Iscr(-m)$ on $\P2q$ with the property  that the 
$\Mscr_i/\Mscr_{i+1}$ are shifted 
line objects. 
\end{theorems} 
\begin{proof} This follows easily from the lemma \ref{ref-5.6.5-70} and 
Corollary \ref{ref-3.6-34}. 
\end{proof} 
\begin{remarks} There is some freedom in choosing  the line objects 
occurring in Theorem \ref{ref-5.6.6-72}. We may assume for example 
that they all map to the same point object. 
\end{remarks} 
\appendix 
\section{Serre duality for graded rings} 
\label{ref-A-74} 
In this section we prove that (a generalization of) BKS-duality holds 
for graded rings. For the convenience of the reader we restate some 
definitions so that this appendix can be read independently of the 
rest of this paper. 

Let $\Ascr$ be a $k$-linear $\Ext$ 
finite triangulated category. By this we mean that for all 
$\Mscr,\Nscr\in\Ascr$ we have $\sum_n\dim_k \Hom(\Mscr,\Nscr[n])<\infty$. 
The category $\Ascr$ is said to satisfy 
Bondal-Kapranov-Serre (BKS) duality  if there is an autoequivalence 
$F:\Ascr\r \Ascr$ together with for all $A,B\in \Ascr$ natural 
isomorphisms 
\[ 
\Hom(A,B)\r \Hom(B,FA)' 
\] 
(where $(-)'$ denotes the $k$-dual). 

Let $\Cscr$ be an  abelian category.  An object $O$ in $D^b(\Cscr)$ is 
said to have finite projective (injective) dimension if 
$\Ext^i(O,\Cscr)=0$ ($\Ext^i(\Cscr,O)=0$) 
for $|i|> u$  for some $u\ge 0$. The minimal such $u$ 
we call the 
projective (injective) dimension of $O$. 

In this appendix we assume that $A$ is a connected graded noetherian 
ring over a  $k$. By $(-)'$ we denote the functor on graded 
vectorspaces which sends $M$ to $\oplus_n M_{-n}^\ast$. If we use 
notations which refer to the left structure of $A$ then we adorn them 
with a superscript ``$\circ$''. 

We make the following additional assumptions on $A$: 
\begin{enumerate} 
\item $A$ satisfies $\chi$ and the functor $\tau$ has finite 
cohomogical dimension. 
\item $A$ satisfies $\chi^\circ$ and the functor $\tau^\circ$ has finite 
cohomogical dimension. 
\end{enumerate} 
These conditions  imply that $A$ has a \emph{balanced dualizing complex} 
\cite{Ye} given by 
$R=R\tau(A)'=R\tau^\circ(A)'$ \cite{VdB16,Ye}. Below we freely use the 
properties of such dualizing complexes. 

We let $D(A)$  be the 
derived category of graded right $A$-modules. $D^b_f(A)$ will be the 
full subcategory of objects in $D^b(A)$ with finitely generated 
homology. The category $D^b_f(\Tails(A))$ is the full subcategory of 
$D^b(\Tails(A))$ consisting of complexes with homology in $\tails(A)$.

We let $D^b_f(A)_{\fpd}$ ($D^b_f(A)_{\fid}$) be the 
full category of $D^b_f(A)$  consisting of objects of finite 
projective (injective) dimension. The categories $D^b_f(\Tails(A))_{\fpd}$  and 
$D^b_f(\Tails(A))_{\fid}$ are defined in a similar 
way. The fact that $\tau$ has finite cohomological dimension implies 
$\pi A(n)\in D^b_f(\Tails(A))_{\fpd}$. 

We will denote the functors $R\underline{\Hom}_A(-,R)$ and 
$R\underline{\Hom}_{A^\circ}(-,R)$ by 
$D$.  Since they define  a duality 
between $D^b_f(A)$ and $D^b_f(A^\circ)$ it is clear that they define 
a duality between 
$D^b_f(A)_{\fid}$ and $D^b_f(A^\circ)_{\fpd}$ and between 
$D^b_f(A)_{\fpd}$ and $D^b_f(A^\circ)_{\fid}$.

It is also clear that these functors induce a duality between 
$D^b_f(\Tails(A))$ and $D^b_f(\Tails(A^\opp))$. We denote these induced 
functors also by  $D$.  Again they define a duality between 
$D^b_f(\Tails(A))_{\fid}$ and $D^b_f(\Tails(A^\circ))_{\fpd}$ and between 
$D^b_f(\Tails(A))_{\fpd}$ and $D^b_f(\Tails(A^\circ))_{\fid}$. 
Recall the following: 
\begin{lemma} \label{ref-A.1-75} 
Let $\Pscr\in  D^b_f(\Tails(A))_{\fpd}$. Then there 
exists an object $P\in D^b_f(A)_{\fpd}$ such that $\Pscr$ is a 
direct summand of $\pi P$. 
\end{lemma} 
\begin{proof} This can be deduced from general results about compact 
objects in triangulated categories. For simplicity we give a direct 
proof based on a trick which the authors learned from Maxim 
Kontsevich. Take $M$ arbitrary such that $\pi M=\Pscr$. 

Take a quasi-isomorphism $Q\r M$ where $Q$ is a 
right bounded complex of finitely generated projective modules. This 
yields a triangle: 
\[ 
(\pi Z)[a] \r\sigma_{\ge -a} \pi Q \r  \Pscr 
\] 
where $Z=\ker (Q_{-a}\r Q_{-a+1})$. 
This triangle corresponds to an element of $\Ext^{a+1}(\Pscr,\pi Z)$ 
which must be zero for large $a$. Hence $\sigma_{\ge a} \pi 
Q=\Pscr \oplus  (\pi Z)[a]$. This proves the lemma. 
\end{proof} 

We recall the following fact. 
\begin{proposition} The functors $-\Lotimes_A R$ and $R\underline{\Hom}_A(R,-)$ 
induce inverse equivalences between $D^b_f(A)_{\fpd}$ and $D^b_f(A)_{\fid}$. 
\end{proposition} 
\begin{proof} If $P\in D^b_f(A)_{\fpd}$ then it is quasi-isomorphic to a 
bounded complex of finitely generated projective $A$-modules. For 
such such a complex 
it is clear that $P\otimes_A R$ has finite injective 
dimension. There is a canonical map $P\r 
R\underline{\Hom}(R,P\otimes_A R)$ which is an isomorphism for 
$P=A$. By induction over triangles one shows that it is an 
isomorphism for all $P$. 

Conversely assume $I\in D^b_f(A)_{\fid}$. Then by duality 
$R\underline{\Hom}(R,I)=R\underline{\Hom}(DI,A)$. 
By the above discussion $DI\in D^b_f(A^\circ)_{\fpd}$. Hence $R\underline{\Hom}_A(DI,A)\in 
D^b_f(A)_{\fpd}$. We 
also find $R\underline{\Hom}_A(DI,A)\otimes_A R=R\underline{\Hom}_A(DI,R)=I$. 
\end{proof} 
The functor $-\otimes_A R$ induces a functor 
$D^-(\Tails(A))\r D^-(\Tails(A)) $ which we denote by $-\otimes 
\Rscr$. Similarly the functor $\RHom_A(R,-)$ induces a functor 
$D^+(\Tails(A))\r D^+(\Tails(A)) $  which we denote by 
$\uRHom(\Rscr,-)$. 

\begin{proposition} 
\label{ref-A.3-76} 
The functors $-\otimes \Rscr$ and $\uRHom(\Rscr,-)$ 
induces inverse equivalences between $D^b_f(\Tails(A))_{\fpd}$ and 
$D^b_f(\Tails(A))_{\fid}$. 
\end{proposition} 
\begin{proof} If $\Pscr\in D^b_f(\Tails(A))_{\fpd}$ then by lemma 
\ref{ref-A.1-75} $\Pscr$ is direct summand of some $\pi P$ with 
$P\in  D^b_f(A)_{\fpd}$. 

Using the proof of the previous proposition 
this easily implies that $\Pscr\otimes \Rscr\in 
D^b_f(\Tails(A))_{\fid}$ and $\uRHom(\Rscr,\Pscr\otimes 
\Rscr)=\Pscr$ (essentially because we may reduce to $\Pscr=\pi A(n)$ 
for some $n$). 

Conversely assume $\Iscr=\pi I \in D^b_f(\Tails(A))_{\fid}$. Then 
$\uRHom(\Rscr,\Iscr)= \linebreak
\pi R\underline{\Hom}(R,I)=\pi R\underline{\Hom}(DI,A)$. 

We have by definition  $\pi DI= D\pi I$, and hence $\pi DI\in 
D^b_f(A)_{\fpd}$. Then it follows from lemma 
\ref{ref-A.1-75} that  $\pi DI$ is a direct summand of some $\pi Q$ 
with $Q\in D_f^b(Q)_{\fpd}$. 

We easily deduce from this that 
$\pi R\underline{\Hom}(DI,A)$ is a direct summand of \linebreak
$\pi R\underline{\Hom}(Q,A)$ and hence $\uRHom(\Rscr,\Iscr)=\pi 
R\underline{\Hom}(DI,A)\in 
D^b_f(\Tails(A))_{\fpd}$. 
The proof now continous as the proof of the previous proposition. 
\end{proof} 

\begin{theorem} (Serre duality) 
For all $\Mscr\in D^b_f(\Tails(A))_{\fpd}$, $\Nscr\in D^b_f(\Tails(A))$ 
there are natural isomorphisms 
\[ 
\Hom(\Mscr,\Nscr)\cong\Hom(\Nscr,F\Mscr)' 
\] 
where 
\begin{equation} 
\label{ref-A.1-77} 
F\Mscr=(\Mscr\Lotimes \Rscr)[-1] 
\end{equation} 
Furthermore the functor $F$ defines an equivalence between 
$D^b_f(\Tails(A))_{\fpd}$ and $D^b_f(\Tails(A))_{\fid}$. 
\end{theorem} 
\begin{proof} 
As in \cite{YZ} our proof of Serre duality is based on the local 
duality formula \cite{VdB16,Ye}. The formulation of local duality in 
\cite{VdB16} 
used the functor $R\tau$ but the same proof works for the functor 
$RQ$ where $Q=\omega\circ \pi$. Furthermore it is possible to throw 
an extra perfect complex 
into the bargain. If we do this we obtain canonical isomorphisms 
\begin{equation} 
\label{ref-A.2-78} 
\underline{\Hom}_A(N,P\otimes_A (RQ A)')\cong \underline{\Hom}_A(P, RQ N)' 
\end{equation} 
for $N\in D(A)$ and $P\in D^b_f(A)_{\fpd}$. By adjointness 
$\underline{\Hom}_A (P, RQ N)_0 = \linebreak
\Hom_{\Tails(A)}(\pi P, \pi N)$. In addition, if we apply 
\eqref{ref-A.2-78} with $N$ finite dimensional then we find 
$\underline{\Hom}_A(N,P\otimes_A (RQ A)')=0$. Thus using lemma 
\ref{ref-A.1-75} we 
obtain for $N\in D^b_f(A)$: $\underline{\Hom}_A(N,P\otimes_A (RQ 
A)')_0=\Hom_{\Tails(A)}(\pi N,\pi (P\otimes_A (RQ 
A)')$. Now the standard triangle for local cohomology yields $RQ 
A=\cone(R\tau A\r A)$ and thus $(RQ A)'=\cone(A'\r  R)[-1]$. Using the 
fact that $A'$ is torsion we easily obtain from this: $\pi (P\otimes_A (RQ 
A)' )=F(\pi P)$ where $F$ is defined as in the statement of the theorem. So now 
we have shown 
\begin{equation} 
\label{ref-A.3-79} 
\Hom_{\Tails(A)}(\pi N,F(\pi P))\cong \Hom_{\Tails(A)}(\pi P, \pi N)' 
\end{equation} 
Now we obtain from 
lemma \ref{ref-A.1-75} that $\Mscr$ 
is a direct summand of a complex $\pi P$ with $P\in D^b_f(A)_{\fpd}$. Thus 
\eqref{ref-A.3-79} is true for $\Mscr$ and this 
finishes of the the first part of the theorem. The last part is 
Proposition \ref{ref-A.3-76}. 
\end{proof} 
\begin{corollary} If $\Tails(A)$ has finite global dimension then 
$D^b_f(\Tails(A))$ satisfies BKS-duality. 
\end{corollary}

\ifx\undefined\bysame 
\newcommand{\bysame}{\leavevmode\hbox to3em{\hrulefill}\,} 
\fi

\end{document}